\documentclass{article}

\usepackage{amsmath,fullpage}
\usepackage{amsfonts}
\usepackage{amsthm}
\usepackage{nicefrac}
\usepackage{dsfont}
\usepackage{amssymb}
\usepackage{mathrsfs}
 \usepackage{graphicx}
\usepackage{xcolor, caption, hyperref, bm, psfrag}

\newtheorem{lm}{Lemma}
\newtheorem{te}{Theorem}
\newtheorem{cor}{Corollary}
\newtheorem{re}{Remark}
\newtheorem{prop}{Proposition}
\newtheorem{con}{Conjecture}

\def\bl{\begin{lm}}
\def\el{\end{lm}}
\def\bprop{\begin{prop}}
\def\eprop{\end{prop}}
\def\bt{\begin{te}}
\def\et{\end{te}}
\def\bc{\begin{cor}}
\def\ec{\end{cor}}
\def\bcon{\begin{con}}
\def\econ{\end{con}}
\def\br{\begin{re}}
\def\er{\end{re}}
\def\be{\begin{equation}}
\def\ee{\end{equation}}
\def\bp{\begin{proof}}
\def\ep{\end{proof}}

\def\RR{{\mathbb R}}
\def\EE{{\mathbb E}}
\def\PP{{\mathbb P}}

\def\VV{{\mathbb V}}
\def\NN{{\mathbb N}}
\def\ML{{\mathcal L}}

\def\MP{{\mathcal P}}
\def\COV{{\mathbb {COV}}}

\def\a{{\alpha}}
\def\b{{\beta}}

\def\1{{\mathbf 1}}

\allowdisplaybreaks

\begin{document}
\sloppy

\thispagestyle{empty}
\title{The total external length of the evolving Kingman coalescent\footnote{Work partially supported by the DFG Priority Programme SPP 1590 ``Probabilistic Structures in Evolution''.  I. D. was partially supported by the German Academic Exchange Service (DAAD).}}
\author{Iulia Dahmer,\thanks{Institut f\"ur Mathematik, Goethe-Universit\"at, 60054 Frankfurt am Main, Germany. \newline\texttt{dahmer@math.uni-frankfurt.de, kersting@math.uni-frankfurt.de}
} \qquad
 G\"otz Kersting$^\dagger$}
\maketitle

\begin{abstract}
The evolving Kingman coalescent is the tree-valued process which records the time evolution undergone by the genealogies of Moran populations. We consider the associated process of total external tree length of the evolving Kingman coalescent and its asymptotic behaviour when the number of leaves of the tree tends to infinity. We show that on the time-scale of the Moran model slowed down by a factor equal to the population size, the (centred and rescaled) external length process converges to a stationary Gaussian process with almost surely continuous paths and covariance function $c(s,t)=\Big( \frac 2 {2+|s-t|} \Big)^2$. A key role in the evolution of the external length is played by the internal lengths of finite orders in the coalescent at a fixed time which behave asymptotically in a multivariate Gaussian manner (see Dahmer and Kersting (2015)). A coupling of the Moran model with a critical branching process is used. We also derive a central limit result for normally distributed sums endowed with independent random coefficients.
\end{abstract}
\begin{small}
\end{small}

\section{Introduction and main result}\label{intro}

In mathematical population genetics the Kingman coalescent is a classical model for describing the genealogies of populations. If their population size is equal to $n$, the $n$-Kingman coalescent can be graphically represented as a binary tree which starts with $n$ leaves and spends an exponential time $X_k$ with parameter $\binom k 2$ having $k$ branches. The inter-coalescence times are independent. When labelling the leaves of the tree by $1, \dots, n$ one can define the $n$-Kingman coalescent as a partition-valued process $\Pi=(\Pi_k)_{1 \leq k\leq n}$ started in the partition $\pi_n=\{\{1\}, \dots, \{n\}\}$ of $\{1, \dots, n\}$ into singletons with the property that, when it is in a state $\pi_k$, it jumps after the exponential time $X_k$ to a state obtained by merging two randomly chosen blocks from $\pi_k$. This happens in the coalescent time direction, from the present to the past.

As time runs forwards the population evolves and its genealogy changes, giving rise to a tree-valued process known as the {\it evolving Kingman coalescent} (\cite{PW06},\cite{PWW11}). We consider in this paper populations started at remote past, so to say at time $-\infty$, and driven by the Moran model. This is a stationary, continuous-time evolution model in which each pair of individuals from the population is picked at rate 1, at which moment one of the individuals dies and the other one gives birth to an offspring, see e.g. \cite{Du08}.

The evolving coalescent discloses features of the Kingman coalescent which are less visible in the static model, when one analyses the tree only at fixed times. It is worth noting that some aspects of the tree arise from the recent past and others from the more distant past and that for certain functionals of the tree it may not be clear which one of these contributions dominates. As we shall see, the influence of the recent and the distant past is reflected among others in different time-scales, namely the evolutionary time-scale and the generations time-scale. We come back to this in more detail below.

\begin{figure}
\psfrag{a}{$t_1$}
\psfrag{b}{$t_2$}
\includegraphics[scale=0.5]{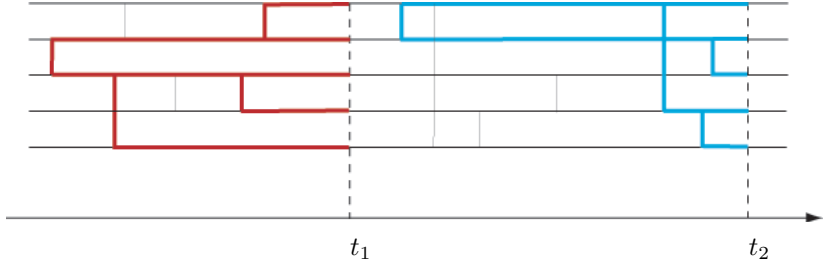}
\caption{A realisation of the evolving Kingman coalescent at two different times for a Moran model with population size $n=5$}
\end{figure}

Particular functionals of coalescent trees such as the total branch length (the sum of the lengths of all the branches of the tree) and the total external length (the sum of the lengths of the external branches) have been extensively investigated in the literature also due to their statistical relevance in population models including mutations. The mutation events are modelled as points of a Poisson process with constant rate on the branches of the coalescent tree. In the Infinitely Many Sites Model the total tree length controls the total number of mutations seen in the population whereas the total external length controls the number of mutations that affect only single individuals.

The Kingman coalescent is one of a big class of coalescent models which have been studied in the literature. The total length and the external length have been extensively studied for the Beta-coalescents. For the asymptotics in the static case we refer to \cite{dr}, \cite{ikmo}, \cite{mo}, \cite{bebesw,bebesw2}, \cite{bebeli}, \cite{ke}, \cite{deyu} and \cite{DKW14a}. In the dynamic case, the evolving Beta$(2-\alpha, \alpha)$-coalescent was introduced and investigated for $\alpha=1$ (the Bolthausen-Sznitman coalescent) in \cite{Sch12} and for $1<\alpha<2$ in \cite{KSW14}. In the latter case the processes of (centred and rescaled) total tree length (for $1<\alpha<\frac 1 2 (1+\sqrt 5)$) and of total external length (for $1<\alpha<2$)  converge in the sense of the finite dimensional distributions to stationary moving average processes with stable distributions. The time-scale is the generations time-scale for both the total length process and the external length process (which is obtained by slowing down time by a factor of  $n^{\alpha-1}$). As we shall see, this is contrary to the Kingman case. Evolving coalescents have also been investigated from a different point of view, namely as evolving metric spaces in (e.g.) \cite{DGP12}, \cite{GPW09}, \cite{GPW13}, \cite{Gu14}.

In this paper we investigate the asymptotic behaviour of the total external length process of the evolving Kingman-coalescent. Observe that this process is {\it stationary}, a property that it inherits from the stationarity of the Moran model. The next theorem states our main result.

\bt\label{mainresult_fdd}
Let $\ML^{n}_t$ be the external length of the evolving Kingman $n$-coalescent at time $t \in \RR$. Then, as $n \to \infty$
\[\left(\sqrt{\frac{n}{4\log n}}\left(\ML^{n}_{\nicefrac{ t}{n}}-2\right)\right)_{ t \in \RR} \stackrel{f.d.d.}{\longrightarrow} \ML,\]
where the limiting process $\ML=\left(\ML_h\right)_{h \in \RR}$ is stationary, Gaussian, a.s. continuous, with mean 0 and covariance function
\[\mathbb{COV}(\ML_0, \ML_{h})=\left(\frac{2}{2+h}\right)^2,  h\geq 0.\]
\et

This convergence for the one-dimensional case can be read off from Janson and Kersting \cite{JK11}.

Let us contrast this result with the theorem of Pfaffelhuber, Wakolbinger and Weisshaupt \cite{PWW11}, who investigated the process $(L^n_t)_{t \in \RR}$ of total tree length of the evolving Kingman $n$-coalescent:

\bt\label{th_PWW}
There exists a stationary process $L=(L_t)_{t \in \RR}$ with paths in $\mathbb D$, the space of c\`adl\`ag functions equipped with the Skorokhod topology, such that 
\[L^n - 2 \log n\stackrel {d} {\longrightarrow} L \quad \text { as } n \to \infty.\]
$L_t$ has a Gumbel distribution for all $t$. Moreover
\[\frac{1}{t |\log t|}\EE((L_t-L_0)^2)\to 4 \text { as } t \to 0. \]
\et

Theorem \ref{th_PWW} deals with the original time-scale from the Moran model, that is the evolutionary time-scale. On the contrary, the time in Theorem \ref{mainresult_fdd} runs on the evolutionary time-scale {\it slowed down} by a factor of $n$. This is what we called the generations time-scale. Note that on the evolutionary time-scale reproduction events happen in the population with a rate of order $n^2$. Therefore, for large populations, after a time of order $\frac 1 n$ each individual will have taken part in a number of reproduction events of order 1. It is not surprising that the external length has to be considered on this time-scale. It might be less obvious that the evolutionary time-scale is the appropriate choice for the total length. 

We use the notation $t$ for the time on the evolutionary time-scale and
\[h=h(n):=\frac{ t}{n},\]
for the time points of the generations time-scale.

Note also that the limiting process $\ML$ is almost surely continuous, whereas $L$ is almost surely made up of jumps. This reflects the fact that the total length experiences big jumps (at the times when old families become extinct). For the external length such extremal events do not come into play. This is also reflected in the type of limiting distributions. Let us point out that the process $L$ is not a semimartingale (\cite{DKW14b}), which so far is an open problem for the process $\ML$.

\vspace{0.5cm}
\noindent
In the rest of this introduction we outline the proof of our theorem. Due to stationarity it suffices to analyse the dynamics of the external length on a time interval $(0,h)$ with $h >0$. 

For describing the dynamics of the total external length of the Kingman coalescent we recall the notion of branches of order $i$ introduced in \cite{DK13}. In the coalescent tree each branch is situated above a subtree. If this subtree has $i$ leaves, one says that the branch is {\it of order i}. If $i=1$ the branches are external and if $i \geq 2$ they are internal. We denote by $\ML^{n,i}_t$ the total length of order $i$ (the sum of the lengths of branches of order $i$) in the $n$-coalescent belonging to the population alive at time $t$.

\begin{figure}
\psfrag{a}{a)}
\psfrag{b}{b)}
\psfrag{c}{$h_1$}
\psfrag{d}{$\zeta$}
\psfrag{e}{$h_2$}
\includegraphics[scale=0.4]{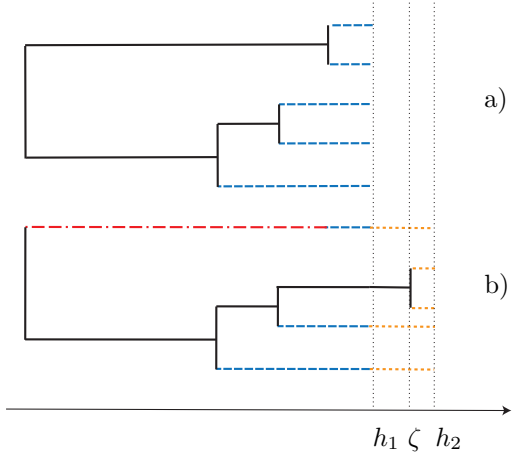}
\caption {In Figure a) the external branches of the coalescent at time $h_1$ are marked in (dashed) blue. At time $\zeta$ a reproduction event happens and the tree changes. The external branches at time $h_2$ are shown in Figure b) with the following code: the (dashed) blue pieces make up the part of the external length from time $h_1$ that is still external at time $h_2$, the red (dash and dot) piece is a branch that was internal at time $h_1$ and became free at time $\zeta$, while the orange (dotted) pieces represent the external length $I^n_{h_1,h_2}$ gathered in the time interval $(h_1,h_2)$}\label{Fig_ext_len_jump}
\end{figure}

Let us now look at the dynamics of the total external length. On the generation time-scale each pair of individuals takes part in a reproduction event at rate $\frac{1}{n}$. As long as no reproduction event takes place, the external branches of the coalescent tree grow linearly and hence the external length grows at rate $1$.  At the time $\zeta$ of a reproducing event, the total external length process has a jump. The length of the branch $b_1$ corresponding to the individual that dies and the length of the branch $b_2$ corresponding to the individual that reproduces are subtracted from the total external length. The branch $b_1$ is removed from the tree, whereas $b_2$ becomes an internal branch of order 2.  Typically, by the removal of $b_1$, a branch that was internal (of order 2) at time $\zeta-$ becomes part of an external branch at time $\zeta$ (the exceptional case is when $b_1$ and $b_2$ stem from the same branching event). The length of this branch is then added to the total external length. Note that this piece may become internal again later on due to new reproducing events. If such an internal branch is part of an external branch at a future time $h$ we say that the branch is {\it free} at time $h$. An example is given in Figure \ref{Fig_ext_len_jump}.

From the tree at time $0$ not only internal branches of order 2 may be free at time $h$, but also internal branches of orders greater than 2. The lengths of the free branches become part of the external length at time $h$ and hence $\ML^n_h=\ML^{n,1}_{h}$ can be written as
\be\label{ext_len_decomp}
\ML^{n}_{h}= \sum_{i=1}^{n-1} \Lambda^{n,i} \ML^{n,i}_{0} + I^n_{0,h},
\ee
where by $\Lambda^{n,i}=\Lambda^{n,i}_{0, h}$ we denote the (random) proportion of the length $\ML^{n,i}_{0}$ of order $i$ that is free by time $h$ and by $I^n_{0,h}$ the external length of the coalescent at time $h$ arising in the time interval $(0,h)$.
 
 \newpage
A key role in the analysis of the external length process is played by the main result of \cite{DK13} which states that for any $h \in \RR$ and any $r \in \mathbb{N}$, as $n \to \infty$
\[
\sqrt{\frac{n}{4 \log n}}\left(\ML^{n,1}_h-\eta_1, \dots, \ML^{n,r}_h-\eta_r\right)\stackrel{d}{\longrightarrow} N\left(0,I_r\right),
\]
where $I_r$ denotes the $r \times r$-identity matrix and 
\be\label{expectation_L}
\eta_i=\EE(\ML^{n,i}_h)=\frac{2}{i} \quad \text{ for every $i \geq 1$.}
\ee

Let us label the individuals alive at time $0$ by $1, \dots, n$ and denote by $\MP_n$ the power set of $\{1, \dots, n\}$. For $A \in \MP_n \setminus \{\emptyset, \{1, \dots, n\}\}$ denote
\begin{align*}
L_A:= &\text{the length of the branch supporting} \text{ the leaves with labels in $A$ in the coalescent at time 0}, \notag
\end{align*}
with the convention $L_A=0$ if there is no such branch in the coalescent and
\[
\mathscr F_A:= \{\text{the number of descendants at time $h$ of the individuals with labels in $A$ is equal to 1}\}.
\]
Observe that with these notations 
\[\1_{\mathscr F_A}\cdot L_A = \1_{\{\text{there is a branch supporting the leaves with labels in $A$ at time 0 which is {\it free} at time } h\}} \cdot L_A\]
and thus the $i$-th summand on the right-hand side of (\ref{ext_len_decomp}) can be rewritten as
\[\Lambda^{n,i} \ML^{n,i}_{0}=\sum_{A \in \MP_n \atop |A|=i} \1_{\mathscr F_A} \cdot L_A.
\]
For any fixed $r \in \NN$ we can write \eqref{ext_len_decomp} as
\be\label{ext_len_decomp_1} 
\ML^{n}_{h}= \sum_{A \in \MP_n \atop 1 \leq |A| \leq r} \1_{\mathscr F_A} \cdot L_A + R^{n,r}_h+ I^n_{0,h},
\ee
where
\be\label{def_R}
R^{n,r}_h:= \sum_{A \in \MP_n \atop r<|A|<n} \1_{\mathscr F_A } \cdot L_A
\ee
records the contribution of the lengths of the branches of orders larger than $r$ to the external length at time $h$.

The proof of Theorem \ref{mainresult_fdd} will show that for large $r$ it is only the first term on the right-hand side of (\ref{ext_len_decomp_1}) that makes a non-negligible contribution to the external length at time $h$. Note moreover that this sum is obtained by a random thinning (by means of the indicators $ \1_{\mathscr F_A}$) of the lengths $\ML^{n,i}_0$ of order at most $r$. In Section \ref{normality_thinned_sum} we prove a general result, which assures that if a sum of (possibly dependent) random variables converging in distribution to a normal-distributed random variable, is endowed with independent random coefficients, then the sum obtained through this procedure continues to be asymptotically Gaussian. The result of \cite{DK13} provides the asymptotic normality assumption on the initial sums $\ML^{n,i}_0$. 

The second key idea of the paper is that the {\it freeing}-mechanism described above can be coupled in a natural way with a birth and death process started with $|A|$ individuals at time 0. This allows to replace the original thinning procedure by another one which provides independent coefficients. Moreover, the coupling is also used in the investigation of the contributions that the branches carrying a large number of leaves make to the evolving external length (in Section \ref{Proof of the Proposition}). The coupling is presented in Section \ref{A birth and death process}.

The third building block of the proof of our main result is a new description of the coalescent tree structure (also called the tree topology). This is obtained by neglecting the exponential lengths of the branches. Then the whole tree structure is captured by the random variables $K_A$ and $J_A$, $A \in \MP_n \setminus \{\emptyset, \{1, \dots, n\}\}$, where $K_A$ is the level in the coalescent at which the branch carrying the leaves with labels in the set $A$ is formed and $J_A$ is the level at which this branch ends (with the convention that $K_A=J_A=\infty$ if no such branch exists in the coalescent). We show that these quantities (and not the exponential times) play the central role in the development of the external length. We examine them in detail in Section \ref{Lower and upper levels}.

Building on this preparatory work and on Fu's estimates \cite{Fu95} we deal in Section \ref{Proof of the Proposition} with the expectation and variance of the length contributed to the evolving external length by the large families from the coalescent at time 0.

The proof of Theorem \ref{mainresult_fdd} is then completed in Section \ref{Proof of the main result}.  Putting the building blocks described above together one more idea is needed for the proof, namely it is necessary to distinguish between the contributions to the evolving external length coming from the possible freeing of the internal lengths that are in small, medium-sized and respectively large depth of the tree. It turns out that for the contributions of the first and the latter groups of branches one can easily replace the freeing procedures by procedures dictated by independent birth and death processes (for each branch separately). In the case of the second group of internal branches (defined as the branches crossing a given fixed level in the tree), it holds that the sets of leaves carried by these branches are all disjoint. This property allows us to apply the coupling from Section \ref{A birth and death process} for all branches in the group simultaneously in order to obtain independent coefficients for the thinning. By the replacement of the freeing-mechanism through independent birth and death processes we also obtain the formula for the covariance of the limiting Gaussian process: $\Big( \frac 2 {2+h}\Big)^2$ is simply the probability that a critical birth and death process starting with one individual at time 0 has one individual at time $h$.

{\it Notation.} We recall the Vinogradov notation:
\begin{align*}
a_n \ll b_n &\text{ if there exists a finite constant $c$ independent of $n$ }
\\
&\text{ (but possibly dependent of $r$ and $h$) such that } |a_n| \leq c b_n.
\end{align*}
We use this notation instead of the O-notation since for voluminous formulas it is easier to read.

Moreover we denote throughout by $c$ a finite constant whose value is not important and
may change from line to line.

\section{Asymptotic normality of sums with random coefficients}\label{normality_thinned_sum}

The following Proposition contains a general statement of the form: If the sum of not necessarily independent random variables is asymptotically normal then this property persists if the summands are supplied with independent random coefficients.
\begin{prop}\label{lemma_ind_thin}
Let $r \in \NN$ be fixed and let $k_1(n), \dots, k_r(n) \in \NN$ be such that $\lim _{n \to \infty} k_i(n) =\infty$ for each $1 \leq i\leq r$. Let $Y_{i,j}, 1 \leq i\leq r, 1 \leq j\leq k_i(n)$ be random variables such that in $\RR^r$
\[ \Big( \sum_{j=1}^{k_i(n)} Y_{i,j} \Big)_{i=1, \dots, r} \stackrel d \longrightarrow N(0,I_r).\]
Also let $U_{i,j}, 1 \leq i\leq r, 1 \leq j\leq k_i(n)$ be independent of one another and independent of the collection $\{Y_{i,j}\}_{1 \leq i\leq r \atop 1 \leq j\leq k_i(n)}$, with the property that for fixed index $i$ the variables  $U_{i,j}, 1 \leq j\leq k_i(n)$ are identically distributed with second moment $m_i$. If the following two conditions hold: for each $1 \leq i\leq r$
\[  \sum_{j=1}^{k_i(n)} Y_{i,j}^2 \stackrel {\PP} \longrightarrow 1 \quad \text{ and } \quad \max_{1 \leq j\leq k_i(n)} |Y_{i,j}| \stackrel {\PP} \longrightarrow 0 \]  then 
\[\sum_{i=1}^r \sum_{j=1}^{k_i(n)} U_{i,j} Y_{i,j} \stackrel d \longrightarrow  N\Big(0, \sum_{i=1}^r m_i\Big).\]
\end{prop}

\bp
The characteristic function of the sum $ \sum_{i=1}^r \sum_{j=1}^{k_i(n)} U_{i,j} Y_{i,j}$ is
\begin{align}\label{char_1}
\varphi^n(\lambda)&=\EE\Big(\EE\Big(\exp\Big(i \lambda  \sum_{i=1}^r \sum_{j=1}^{k_i(n)} U_{i,j} Y_{i,j}\Big) \mid Y_{i,j}, i\leq r, j\leq k_i(n)\Big)\Big) =\EE\Big(\prod_{i=1}^r\prod_{j=1}^{k_i(n)}\varphi_{U_{i,j}}(\lambda Y_{i,j})\Big). 
\end{align}

Using the fact that for complex numbers $\{a_k\}_{1\leq k\leq n}$ and $\{b_k\}_{1\leq k\leq n}$ with $|a_k|\leq 1$ and $|b_k|\leq 1$ for all $k$ it holds that 
\[\Big | \prod_{k=1}^n a_k - \prod_{k=1}^n b_k\Big| \leq \sum_{k=1}^n |a_k-b_k|\]
we obtain by denoting the expectation and variance of $U_{i,j}$ by $\mu_i$ and $\sigma_i^2$ respectively that
\begin{align*}
\Big|\prod_{i,j} \varphi_{U_{i,j}}(\lambda Y_{i,j}) -& \prod_{i,j} \exp\Big(i \mu_i\lambda Y_{i,j} -\frac{\sigma^2_i}{2}\lambda^2Y_{i,j}^2 \Big) \Big|\leq\sum_{i,j}\Big|\varphi_{U_{i,j}}(\lambda Y_{i,j})- \exp\Big(i \mu_i\lambda Y_{i,j} -\frac{\sigma^2_i}{2}\lambda^2Y_{i,j}^2 \Big) \Big| 
\end{align*}
and by using a Taylor expansion that
\begin{align*}
\Big|\prod_{i,j} \varphi_{U_{i,j}}(\lambda Y_{i,j}) - \prod_{i,j} \exp\Big(i \mu_i\lambda Y_{i,j} -\frac{\sigma^2_i}{2}\lambda^2Y_{i,j}^2 \Big) \Big|& \leq  \sum_{i,j}\Big| 1+i \mu_i \lambda Y_{i,j} -\frac{1 }{2}m_i\lambda^2Y_{i,j}^2 + o(|Y_{i,j}|^2) 
\\
&\hphantom{12}-\Big(1+i\mu_i\lambda Y_{i,j} -\frac{\sigma_i^2}{2}\lambda^2 Y_{i,j}^2 -\frac{\mu_i^2}{2}\lambda^2Y_{i,j}^2 + O(|Y_{i,j}|^3) \Big)\Big| \notag
\\
& =\sum_{i,j} o(|Y_{i,j}|^2).
\end{align*}

By assumption $\sum_{i,j} o(|Y_{i,j}|^2)$ converges in probability to 0 as $n \to \infty$. It follows by dominated convergence that
\begin{align}\label{char_2}
\EE\Big(&\Big|\prod_{i,j} \varphi_{U_{i,j}}(\lambda Y_{i,j}) - \prod_{i,j} \exp\Big(i \mu_i\lambda Y_{i,j} -\frac{\sigma^2_i}{2}\lambda^2Y_{i,j}^2 \Big) \Big|\Big) \stackrel {n \to \infty} \longrightarrow  0.
\end{align}
From the assumptions on the collection of random variables $\{Y_{i,j}\}_{1 \leq i\leq r \atop 1 \leq j\leq k_i(n)}$ it follows that $\Big(\sum_j Y_{i,j}, \sum_j Y^2_{i,j} \Big)_i$ converges in distribution and we obtain from \eqref{char_1} and \eqref{char_2}  that 
\begin{align*}
\lim_{n \to \infty} \varphi^n(\lambda) &= \lim_{n \to \infty}\EE\Big(\prod_{i=1}^r\prod_{j=1}^{k_i(n)} \exp(i \mu_i\lambda Y_{i,j} -\frac{\sigma^2_i}{2}\lambda^2Y_{i,j}^2 ) \Big)
\\
&=\lim_{n \to \infty}\EE\Big( \exp(i \lambda \sum_{i=1}^r \mu_i \sum_{j=1}^{k_i(n)} Y_{i,j} - \sum_{i=1}^r \lambda^2 \frac{\sigma_i^2}{2} \sum_{j=1}^{k_i(n)} Y_{i,j}^2) \Big) 
\\
& = e^{-\frac{\lambda^2}{2} \sum_{i=1}^r \mu_i^2} \cdot e^{-\frac {\lambda^2 } 2 \sum_{i=1}^r \sigma_i^2} 
\\
&= e^{-\frac{\lambda^2}{2} \sum_{i=1}^r m_i}.
\end{align*}
Since the right hand side is the characteristic function of the normal distribution with expectation 0 and variance $ \sum_{i=1}^r m_i$ the proof of the proposition is finished.
\ep

\section{Coupling a critical branching process to the Moran model}\label{A birth and death process}

In the time interval $(0,h)$ on the generations time-scale the individuals in the population are subject to binary reproduction events and death events according to the Moran dynamics with rate $\frac 1 n$ per pair. In this section we give the description of the evolution of the number of descendants of the individuals alive at time 0 in terms of independent linear birth and death processes. 

We call the descendants of the individual with label $j \in \{1, \dots, n\}$ at time 0, the {\it $\{j\}$-family} and we are interested in the size of the family (number of individuals at each time point). We denote by
\[Z^j:=(Z^j_h)_{h \geq 0}\]
the process recording the size of the $\{j\}$-family as time evolves. Then $Z^j_0=1$ and according to the Moran dynamics of the population, the process $Z^j$ is a continuous time Markov process with jump rates
\[r_{k, k'}=\begin{cases} \frac{1}{2}k(1-\frac{k}{n}) & \text{ if } k'=k+1 
					\\   \frac{1}{2}k(1-\frac{k}{n}) & \text{ if } k'=k-1.
		\end{cases}\]
Indeed, it is only the reproduction events that happen between a member of the $\{j\}$-family and a non-$\{j\}$-family member, that affect the size $Z^j$ of the $\{j\}$-family. If $Z^j$ is in state $k$ there are $k(n-k)$ such pairs and the reproduction events arise at rate $\frac 1 n$. The family size may either decrease by one or increase by one at the time of such a reproduction event, depending on whether the individual from the $\{j\}$-family dies or reproduces. The reproduction events between two $\{j\}$-family members or between two non-$\{j\}$-family members leave $Z^j$ unchanged. 

Let us now fix $A \in \MP_n$ and let $i=|A|$, $1 \leq i \leq n-1$. Denote by 
\[Z^A=(Z^A_h)_{h \geq 0}\]
the process which records the size of the $A$-family as time evolves,
\[Z^A:=\sum_{j \in A} Z^j.\]

To the processes $Z^j$, $j \in A$, we now couple jointly linear birth and death processes $B^{j}=B^{j,A}$. To this aim, let us first observe that the reproduction events that affect the size $Z^j$ of the $\{j\}$-family, $j \in A$, can be grouped as follows: on one hand there are the events between an individual from the $\{j\}$-family and an individual from the $\big(\{1, \dots, n\} \setminus A\big)$-family (we call these events {\it $\a$-events}) and on the other hand there are the reproduction events between an individual from the $\{j\}$-family and an individual from the $\big(A\setminus \{j\}\big)$-family (we call these events {\it $\b$-events}). If the process $Z^j$ is in state $k$ and the process $Z^A$ is in state $l$ then the $\alpha$-events occur at rate  $k(n-l)/n$, whereas the $\beta$-events occur at rate $k(l-k)/n$.

We construct the coupled processes $B^j$ starting from the processes $Z^j$, $j \in A$, in two steps. The first step is to leave aside the influence of the $\beta$-events on the size of the $\{j\}$-family. This gives rise to a process $Y^j$ which, when $Z^A$ is in state $l$, jumps at rate
\[\tilde r_{k, k'}=\begin{cases} \frac{1}{2}k(1-\frac{l}{n}) & \text{ if } k'=k+1 
					\\   \frac{1}{2}k(1-\frac{l}{n}) & \text{ if } k'=k-1.
		\end{cases}\]

The second step of the construction is to perform a random time change $\theta_h=\theta_h^A$ on the process $Y^j$ depending on the state of the process $Z^A$. Note that $n$ is an absorbing state for the process $Z^A$ and denote by
\[\tau=\tau^A:=\inf\{u \geq 0: Z^A_u=n\}\]
the hitting time of the process $Z^A$ at state $n$.

Let
\be\label{def_theta}
\theta_h=\theta_h^A:=\int_0^h \Big(1-\frac{Z^A_u}{n}\Big)du.
\ee
Observe that $\theta_h$ as a function of $h$ is continuous and strictly increasing for $h \leq \tau$. With this notation we define $B^j_{\theta}$ for $\theta\leq \theta_\tau$ implicitly by
\begin{align}\label{def_proc_B}
B^j_{\theta_h}=B^{j,A}_{\theta_h}:=Y^j_h, &\quad  h \leq \tau.
\end{align}
Note that $\theta_h<h$ for all $0<h\leq \tau$, thus $B^j$ is a randomly speeded-up version of the process $Y^j$. The process $B^j$ is (using results of \cite{EK86} Section 6.1) a birth and death process with  birth and death rates both equal to $\frac{1}{2} k$ when the process is in state $k$. In order to fully define the processes $B^j$, $j \in A$, we let them evolve from the states $B^j_{\theta_\tau}$ according to these rates independently (of each other and of everything else).

\vspace{0.5cm}
Observe from the construction given above for the processes $B^j$, $j \in A$, that the $\beta$-events (which are eliminated for each $j$ in order to obtain the process $Y^j$ from $Z^j$) do not affect the size $Z^A$ of the $A$-family since they take place between two individuals of this family (one individual reproduces and another dies and the size of the $A$-family remains constant). Thus
\be\label{Z_A=Y_A} 
Z^A=\sum_{j \in A} Z^j = \sum_{j \in A} Y^j=:Y^{A}
\ee
and as a consequence by \eqref{def_proc_B} for all $h\leq \tau$
\be\label{B_and_Z}
B^A_{\theta_h}=Z^A_h,
\ee
where $B^A:=\sum_{j \in A} B^j.$

Define the $\sigma$-algebra
\[\mathcal F_{h}:= \sigma(B^j_s, s \leq h, j \in A)\]
and observe that
\[\theta_\tau=\inf\{u \geq 0: B^A_u=n\}.\]
Thus, $\theta_\tau$ is a stopping time with respect to the filtration $\mathcal F:=(\mathcal F_h)_{h\geq 0}$. For $h \leq \tau$ the relation \eqref{def_theta} can also be expressed as
\be\label{theta_h_deriv_2}
h = \int_0^{\theta_h} \frac 1 {1 - \frac {B^A_u} n}du
\ee
and thus for any $h \geq 0$ it holds
\[\{\theta_h <s\}=\{\theta_{\tau} <s\} \cup\Big\{\theta_{\tau} \geq s, \int_0^s  \frac 1 {1 - \frac {B^A_u} n}du >h\Big\} \in \mathcal F_s.\]
Therefore $\theta_h$ is also a stopping time with respect to the filtration $\mathcal F$.

\bl\label{indep_B_j}
The processes $B^j$, $j \in A$, defined above are independent birth and death processes and also independent of the coalescent tree at time 0.
\el

\bp
Without loss of generality we consider $A=\{1, \dots, i\}$. Then from the construction given above, since we let the $\beta$-events aside, $Y=(Y^1, \dots, Y^i)$ is a Markov process which jumps, when $Z^A=Y^A$ is in state $l$, at rate
\[\lambda_{(k_1, \dots, k_j, \dots, k_i),(k_1, \dots, k'_j, \dots, k_i)}=\begin{cases} \frac{1}{2}k_j(1-\frac{l}{n}) & \text{ if } k'_j=k_j+1 
					\\   \frac{1}{2}k_j(1-\frac{l}{n}) & \text{ if } k'_j=k_j-1.
		\end{cases}\]
By performing the random time change from above on the process $Y$ one then obtains the Markov process $(B^1, \dots, B^i)$ with jumps rates
\[\tilde \lambda_{(k_1, \dots, k_j, \dots, k_i),(k_1, \dots, k'_j, \dots, k_i)}=\begin{cases} \frac{1}{2}k_j & \text{ if } k'_j=k_j+1 
					\\   \frac{1}{2}k_j& \text{ if } k'_j=k_j-1.
		\end{cases}\]
This implies the first assertion. The independence of the coalescent tree at time 0 follows from the independence properties of the underlying Moran model.
\ep

Before further analysing the relation between the processes $B^j$ and $Z^j$, $j\in A$, we compute a few basic quantities connected to $B^j$, that we will use in the sequel. Let 
\[F(s,h)=\sum_{k=0}^\infty s^k \PP(B^j_h=k)\]
denote the generating function of the process $B^j$. Then (see \cite{AN72}, Chapter III) $F$ is the unique solution of the backward equation 
\[\frac{\partial}{\partial h}F(s,h) = \frac{1+F^2(s,h)}{2}-F(s,h)=\frac 1 2 (1-F(s,h))^2\]
with boundary condition $F(s,0)=s$. Integrating the equation gives
\be\label{gen_function}
F(s,h) = 1- \frac{1-s}{1 +\frac{h}{2}-\frac{h}{2}s},
\ee
which can be rewritten as
\[F(s,h) = 1- \frac{1-s}{1 +\frac{h}{2}}\cdot \sum_{j=0}^\infty \Big(\frac{h}{h+2}\cdot s\Big)^j.\]
This gives 
\be\label{BD_probs}
\PP(B^j_h=0)= \frac{h}{h+2} \quad \text{ and }\quad \PP(B^j_h=1)= \Big(\frac{2}{h+2}\Big)^2.
\ee

The process $B^{A}$ has the generating function $F^i$ and it holds that
\be\label{BD_prob_i}
\PP(B^{A}_h=1)= i\cdot \Big(\frac{h}{h+2}\Big)^{i-1} \Big(\frac{2}{h+2}\Big)^2.
\ee

Observe that for fixed $h$ the generating function $F$ and all its derivatives exist for all $s < \frac {h+2} h$ and therefore all moments of $B^{A}_h$ exist and are finite functions of $h$. For $d \in \NN$
\[\EE\Big((B^{A}_h)^{d}\Big)= \EE\Big(\Big(\sum_{j \in {A}}B^j_h\Big)^{d}\Big).\]
Expanding the sum it follows that
\be\label{polinomial}
\EE\Big((B^{A}_h)^d\Big) \leq c_h \cdot i^d,
\ee
for a constant $c_h$ depending on $h$.

\vspace{1cm}
In the rest of this section we focus on the relationship between the processes $Z^{j}$ and $B^{j}$. For $A' \subset A$, $|A'|=i'$ denote similarly as before $Z^{A'}=\sum_{j \in A'} Z^j$, $Y^{A'}=Y^{A',A}=\sum_{j \in A'} Y^j$ and $B^{A'}=B^{A', A}=\sum_{j \in A'} B^j$. Note that contrarily to $Z^{A'}$, the distributions of the processes $Y^{A'}$ and $B^{A'}$ depend on the set $A$ which we fixed in advance.
\bl\label{lemma_probs}
For $A' \subset A$, $|A'|=i'\leq i$ the following bounds hold
\be\label{prob_Y'}
\PP(Y^{A'}_h=1) \leq 2i'\cdot \Big(\frac{h}{h+2}\Big)^{i'-1}
\ee
and 
\be\label{prob_different}
\PP\big( B^{A'}_h= 1, Z^{A'}_h\neq 1\big)+\PP\big(Z^{A'}_h=1, B^{A'}_h\neq 1\big) \leq c \cdot \frac i {n} \cdot \gamma^{i'}
\ee
for some $\gamma <1$ and $c<\infty$ depending only on $h$.
\el

When $A'=A$ the inequality \eqref{prob_Y'} becomes
\be\label{expectation_1_SA}
\PP(Z^A_h=1) \leq 2i\cdot \Big(\frac{h}{h+2}\Big)^{i-1}.
\ee

\bp
We write for short $Z^{A'}=Z'$, $Y^{A'}=Y'$ and $B^{A'}=B'$ and define the random variables
\begin{align*} N(Y',h):=  &\text{ the number of the individuals in $A'$ that have 0 descendants at time $h$, }
\\
&\qquad \text{ when evolving according to $Y'$} 
\end{align*}
and $N(B', h)$ analogously.

Observe from the definition of $Y^j, j \in A$, that from the time $\tau$ on, the state of the process $Y^j$ remains unchanged, that is
\[Y^j_h=Y^j_\tau \qquad \text{ for all } h \geq \tau.\]
Thus, by \eqref{def_proc_B} it follows that for any $h \geq 0$
\be\label{from_Y_to_B}
Y^j_h= B^j_{\theta_h \wedge \theta_{\tau}}.
\ee
Since $Y'_0=i'$ and $\theta_h \wedge\theta_{\tau} < h$ we obtain that for any $h \geq 0$
\begin{align*}
\{Y'_h=1\} & \subset \{N(Y', h)=i'-1\} = \{N(B', \theta_h \wedge\theta_{\tau})= i'-1\}  \subset \{N(B', h)\geq i'-1\}.
\end{align*}
Using \eqref{BD_probs} we thus obtain that 
\[
\PP(Y'_h=1)\leq i'\cdot \Big(\frac{h}{h+2}\Big)^{i'-1}+\Big(\frac{h}{h+2}\Big)^{i'} \leq 2i'\cdot \Big(\frac{h}{h+2}\Big)^{i'-1},
\]
which gives the first claim of the lemma.

\vspace{0.5cm}

As to the second claim observe that the $\beta$-events between a member of the $\{A'\}$-family and a member of the $(A \setminus A')$-family are the $\beta$-events that lead to changes in the size $Z'$ of the $A'$-family. We call these events $\beta^{A'}${\it-events} or for short $\beta'${\it-events}. These occur at rate $k(l-k)/n$ when the processes $Z'$ and $Z^A$ are in the states $k$ and $l$ respectively. 

Let 
\[
\rho=\rho^{A, A'}:=\inf\{u \geq 0: \text{a $\beta'$-event takes place at time $u$}\}
\]
and
\[
E=E_{\beta', h}:=\{\rho <h \}.
\]
Note that on $E^c$ it holds that 
\be\label{Z_A'_Y_A'_on_complement}
Z'_s=Y'_s \quad\text{ for all } s \in [0,h).
\ee
Then
\begin{align}\label{term_1}
\PP\big( &B'_h= 1, Z'_h\neq 1\big)+\PP\big(Z'_h=1, B'_h\neq 1\big) \notag
\\
&\leq\PP\big( B'_h= 1, E\big)+\PP\big( B'_h= 1, Y'_h \neq 1 \big) + \PP\big( Z'_h= 1, E\big)+\PP\big( Y'_h= 1, B'_h \neq 1 \big). 
\end{align}

Observe that since the $\beta'$-events do not affect the evolution of $Y'$ and $Y^A$, in view of \eqref{def_proc_B}, \eqref{Z_A=Y_A}, \eqref{theta_h_deriv_2} and \eqref{Z_A'_Y_A'_on_complement}, given $Y'$ and $Y^A$, the first $\beta'$-event occurs at rate $\frac {Y'_s (Y^A_s-Y'_s)} n \cdot \exp\Big(-\int_0^s \frac {Y'_u (Y^A_u-Y'_u)} n du\Big) ds$. Thus for the first term on the right-hand side of \eqref{term_1} we obtain
\begin{align}\label{rate}
\PP\big( B'_h= 1, E\big)&=\EE\Big(\PP \big(E \mid \mathcal F_{h}\big); B'_h= 1 \Big) \notag
\\
& \leq \EE\Big(\int_0^h \frac {Y'_s (Y^A_s-Y'_s)} n ds; B'_h= 1 \Big) 
\\
& \leq \frac 1 n \EE\Big(\int_0^h Y'_s Y^A_s ds; B'_h= 1 \Big) \notag
\end{align}
and by Fubini's theorem and H\"older's inequality
\begin{align*}
\PP\big( B'_h= 1, E\big)& \leq \frac 1 n \int_0^h \EE \big((Y'_s)^3\big)^{\nicefrac 1 3} \cdot\EE \big((Y^A_s)^3\big)^{\nicefrac 1 3} \cdot \PP(B'_h= 1)^{\nicefrac 1 3}ds.
\end{align*}

Since the function $x \mapsto x^3$ is convex for $x\geq 0$ and $B'$ is a martingale with respect to the filtration $\mathcal F$, it follows that the process $(B')^3$ is a submartingale. We use \eqref{from_Y_to_B} together with the fact that $\theta_s$ and $\theta_{\tau}$ are stopping times and that $\theta_s \leq s$ and obtain for $s \leq h$ that
\begin{align}\label{prod_5_0}
\EE((Y'_s)^3) = \EE((B'_{\theta_s \wedge \theta_{\tau}})^3) \leq \EE((B'_{\theta_s})^3) \leq \EE((B'_{s})^3)\leq  \EE((B'_{h})^3).
\end{align}
Using \eqref{polinomial} we have that for $s \leq h$
\be\label{prod_5}
\EE( (B'_s)^3) \leq c \cdot {i'}^3 \qquad \text{and} \qquad \EE( (Y'_s)^3) \leq c \cdot {i'}^3,
\ee
in particular
\be\label{prod_2}
\EE( (Y^A_s)^3) \leq c \cdot i^3
\ee
for $c<\infty$ depending on $h$.
Now using \eqref{BD_prob_i} it follows that
\[
\PP\big( B'_h= 1, E\big) \leq c\cdot  \frac {i} {n} \cdot i' \cdot \gamma^{i'}
\]
for some $\gamma<1$ and a finite constant $c$ depending on $h$. By further enlarging $\gamma$ we can drop the factor $i'$ and thus obtain that 
\be\label{term_1_1}
\PP\big( B'_h= 1, E\big) \leq c\cdot  \frac {i} {n}\cdot \gamma^{i'}
\ee
for $\gamma<1$ and $c<\infty$ depending on $h$

For the third term on the right-hand side of \eqref{term_1} observe by using the strong Markov property and the bound for $\PP(Z'_h=1)$ given in \eqref{expectation_1_SA} for $A'=A$ that
\begin{align*}
\PP(Z'_h=1, E \mid Z'_s, s \leq \rho)= \1_E \cdot \PP_{Z'_\rho}(Z'_{h-\rho}=1) \leq \1_E \cdot 2 Z'_\rho \Big(\frac {h-\rho}{h-\rho+2}\Big)^{Z'_{\rho}-1}
\end{align*}
and thus for a finite constant $c$ depending on $h$
\begin{align*}
\PP(Z'_h=1, E)\leq c \cdot \EE\Big(\1_E \cdot Z'_\rho \Big(\frac {h}{h+2}\Big)^{Z'_{\rho}}\Big) = c \cdot \EE\Big(\1_E \cdot  \EE\Big(Z'_\rho \Big(\frac {h}{h+2}\Big)^{Z'_{\rho}} \mid Z'_s, s < \rho \Big)\Big).
\end{align*}
Since $Z'_s=Y'_s$ for all times $s<\rho$ and at the time of a $\beta'$-event the process $Z'$ jumps by either +1 or -1 with equal probability and independently of the past we further deduce that
\begin{align*}
\PP(Z'_h=1, E)&\leq c \cdot \EE\Big(\1_E \cdot \Big(\frac 1 2 (Y'_\rho+1) \Big(\frac {h}{h+2}\Big)^{Y_{\rho}+1} +  \frac 1 2 (Y'_\rho-1) \Big(\frac {h}{h+2}\Big)^{Y'_{\rho}-1} \Big)\Big)
\\
& \leq c \cdot \EE\Big(\1_E \cdot (Y'_\rho+1) \Big(\frac {h}{h+2}\Big)^{Y'_{\rho}}\Big)
\end{align*}
and similarly as in \eqref{rate}
\begin{align*}
\PP(Z'_h=1, E)& \leq c \cdot \EE\Big(\int_0^h (Y'_s+1) \Big(\frac {h}{h+2}\Big)^{Y'_{s}} \cdot \frac 1 n \cdot Y'_s (Y^A_s-Y'_s)ds\Big)
\\
& \leq c \cdot \frac 1 n \cdot \EE\Big(\int_0^h {Y'}^2_s Y^A_s \cdot \Big(\frac {h}{h+2}\Big)^{Y'_{s}} ds\Big).
\end{align*}
Now by Fubini's theorem and H\"older's inequality we obtain that
\begin{align}\label{term_3_intermed}
\PP(Z'_h=1, E)& \leq c \cdot \frac 1 n \int_0^h \EE \big((Y'_s)^6\big)^{\nicefrac 1 3} \cdot\EE \big((Y^A_s)^3\big)^{\nicefrac 1 3} \cdot \EE \Big(\Big(\frac {h}{h+2}\Big)^{3Y'_{s}}\Big)^{\nicefrac 1 3}ds.
\end{align}

Since the function $x \mapsto \Big(\frac h {h+2}\Big)^{3x}$ is convex and the process $B'$ is a martingale, the process $\Big( \frac {h}{h+2}\Big)^{3B'}$ is a submartingale. We thus get by \eqref{from_Y_to_B} that for $s \leq h$
\[\EE \Big(\Big(\frac {h}{h+2}\Big)^{3Y'_{s}}\Big)=\EE\Big(\Big( \frac {h}{h+2}\Big)^{3B'_{\theta_s\wedge \theta_\tau}}\Big)\leq\EE\Big(\Big( \frac {h}{h+2}\Big)^{3B'_{h}}\Big)\Big)\]
and moreover
\be\label{prod_4}
\EE\Big(\Big( \frac {h}{h+2}\Big)^{3Y'_{s}}\Big) \leq F\Big( \Big( \frac {h}{h+2}\Big)^3, h\Big)^{i'}
\ee
where the generating function $F$ is given in \eqref{gen_function}. Now using the fact that $F(u,h)<1$ for $u<1$ we obtain by plugging this together with \eqref{prod_5} and \eqref{prod_2} in \eqref{term_3_intermed} that
\be\label{term_1_3}
\PP\big( Z'_h= 1, E\big) \leq c\cdot  \frac {i} {n} \cdot \gamma^{i'}. 
\ee
for $\gamma <1$ and $c<\infty$ depending on $h$.

Let us now turn to the second term on the right-hand side of \eqref{term_1}. By \eqref{from_Y_to_B} it holds that
\begin{align*}
\PP\big( B'_h= 1, Y'_h \neq 1 \big) & = \PP\big( B'_h= 1, B'_{\theta_h \wedge \theta_\tau} \neq 1\big) = \EE\big(\PP(B'_h=1 \mid \mathcal F_{ \theta_h\wedge \theta_\tau}); B'_{\theta_h \wedge \theta_\tau} \neq 1 \big)
\end{align*}
and by \eqref{BD_prob_i}
\begin{align}\label{intermediate}
\PP\big( B'_h= 1, Y'_h \neq 1 \big) & \leq \EE\Big( B'_{\theta_h \wedge \theta_\tau} \Big( \frac {h-\theta_h\wedge \theta_\tau}{h-\theta_h\wedge \theta_\tau+2}\Big)^{B'_{\theta_h\wedge \theta_\tau}-1}; B'_{\theta_h\wedge \theta_\tau}\geq 2\Big)
\\
& \leq \EE\Big((h-\theta_h\wedge \theta_\tau)\cdot B'_{\theta_h \wedge \theta_\tau}\cdot \Big( \frac {h}{h+2} \Big)^{B'_{\theta_h \wedge \theta_\tau}-2}\Big). \notag
\end{align}
By H\"older's inequality we obtain that
\begin{align}\label{term_2}
\PP\big( B'_h= 1, Y'_h \neq 1 \big) & \leq \EE\big((h-\theta_h\wedge \theta_\tau)^3\big)^{\nicefrac 1 3}\cdot\EE\big(\big(B'_{\theta_h \wedge \theta_\tau}\big)^{3}\big)^{\nicefrac 1 3} \cdot\EE\Big(\Big( \frac {h}{h+2}\Big)^{3(B'_{\theta_h \wedge \theta_\tau}-2)} \Big)^{\nicefrac 1 3}. 
\end{align}

Now using that $Z^A=Y^A$ and the fact that the state of the process $Y^j$, $j \in A$, remains unchanged from time $\tau$ on, we get that
\begin{align*}
h-\theta_h\wedge \theta_\tau &\leq h-\theta_h + h \cdot \1_{\{\theta_h >\theta_\tau\}} \leq h-\theta_h + h \cdot \frac {Y^A_h}{n}.
\end{align*}
By \eqref{def_theta} and \eqref{B_and_Z}
\begin{align*}
h-\theta_h\wedge \theta_\tau &\leq \int_0^h \frac {B^A_{\theta_s}} {n} ds + h \cdot \frac {Y^A_h}{n} = \int_0^h \frac {B^A_{\theta_s}+Y^A_h} {n} ds.
\end{align*}
Using Jensen's inequality and Fubini's theorem it follows that
\begin{align*}
\EE\Big((h-\theta_h\wedge \theta_\tau)^3\Big)& \le \EE\Big(h^3\int_0^{h}\Big(\frac {B^A_{\theta_s}+Y^A_h} {n} \Big)^3\frac {ds} h \Big) =\frac {h^2}{n^3} \int_0^h  \EE\Big( \big(B^A_{\theta_s}+Y^A_h\big)^3 \Big) ds.
\end{align*}
By \eqref{prod_5_0} and \eqref{prod_5} (for $A'=A$) together with \eqref{prod_2} we get that for $c <\infty $ depending on $h$
\begin{align}\label{prod_3_3}
\EE\Big((h-\theta_h\wedge \theta_\tau)^3\Big)& \leq c \cdot \frac {i^3} {n^3}.
\end{align}

By the same arguments used to obtain \eqref{prod_4} and using again \eqref{prod_5_0} and \eqref{prod_5} in \eqref{term_2} we obtain that
\be\label{term_1_2}
\PP\big( B'_h= 1, Y'_h \neq 1 \big) \leq c \cdot \frac 1 {n}\cdot i  \cdot \gamma^{i'} 
\ee
for $\gamma <1$ and $c<\infty$ depending on $h$.

For the last term on the left-hand side of \eqref{term_1} we deduce that 
\begin{align*} 
\PP\big(Y'_h=1, B'_h\neq 1\big)& =  \PP\big(B'_{\theta_h\wedge \theta_\tau }=1, \text{ $\exists$ at least one birth or death event in } (\theta_h\wedge \theta_\tau, h)\big) 
\\
& = \EE\Big(\PP\big(\text{ $\exists$ at least one birth or death event in } (\theta_h\wedge \theta_\tau, h) \mid \mathcal F_{\theta_h\wedge \theta_\tau}\big); B'_{\theta_h\wedge \theta_\tau} =1 \Big)  \notag
\\
& \leq\EE\Big(1-e^{-(h-\theta_h\wedge \theta_\tau)}; B'_{\theta_h\wedge \theta_\tau } =1 \Big) 
\\
& \leq \EE\big(h-\theta_h\wedge \theta_\tau; Y'_h =1).
\end{align*}
Using H\"older's inequality together with \eqref{prod_3_3} and \eqref{prob_Y'} we obtain that
\be\label{term_1_4}
\PP\big(Y'_h=1, B'_h\neq 1\big)\leq c \cdot \frac 1 {n} \cdot i \cdot \gamma^{i'} 
\ee
for some $\gamma <1$ and $c<\infty$ depending on $h$. 

This together with \eqref{term_1_1}, \eqref{term_1_3} and \eqref{term_1_2} gives the second claim of the lemma.
\ep

\vspace{1cm}
Let the set $A$ continue to be fixed, consider $A'\subset A$ with $|A'|=i'$ and denote $A'':=A\setminus A'$ and $i'':=i-i'=|A''|$.  As before consider the associated processes $Z'':=Z^{A''}$, $Y'':=Y^{A''}$ and $B'':=B^{A''}$.

\bl\label{lemma_probs_2}
The following bound holds
\begin{align}\label{bound_four prods}
\PP\big(B^{A'}_h&= 1, Z^{A'}_h\neq 1,  Z^{A''}_h=1\big)+\PP\big(Z^{A'}_h=1, B^{A'}_h\neq 1,  Z^{A''}_h=1\big)
\\
& + \PP\big( B^{A'}_h= 1, Z^{A'}_h\neq 1, B^{A''}_h=1\big)+\PP\big(Z^{A'}_h=1, B^{A'}_h\neq 1,  B^{A''}_h=1\big) \notag
\\
& \qquad \leq c \cdot \frac 1 n \cdot \gamma^{i}. \notag
\end{align}
for some $\gamma<1$ and $c<\infty$ depending on $h$.
\el 

\bp
We use similar arguments as in the proof of Lemma \ref{lemma_probs}. Note from definition of the $\beta^{A'}$-events (which are actually also the $\beta^{A''}$-events) that these events involve an individual from the $A'$- and one from the $A''$-family and that they lead to changes in the sizes of the two families. 

It holds that 
\[\{B'_h=1, Z'_h\neq 1, Z''_h=1 \} \subset  \{B'_h=1, Z''_h=1, E\} \cup \{B'_h=1, Y''_h=1, Y'_h\neq B'_h\}\]
and similarly for the other terms in (\ref{bound_four prods}). Therefore the left-hand side of \eqref{bound_four prods} is bounded from above by
\begin{align}\label{cov_1_SA_second_ineq}
\PP&\Big(B'_h=1, Z''_h=1,E\Big) + \PP\Big(B'_h=1,Y''_h=1, Y'_h\neq B'_h\Big)
\\  
&  \hphantom{12} + \PP\Big(Z'_h=1, Z''_h=1, E\Big) + \PP\Big(Y'_h=1, Y''_h=1, Y'_h\neq B'_h\Big) \notag
\\  
& \hphantom{12} + \PP\Big(B'_h=1, B''_h=1, E\Big) + \PP\Big(B'_h=1, B''_h=1, Y'_h\neq B'_h\Big) \notag
\\  
& \hphantom{12} + \PP\Big(Z'_h=1, B''_h=1, E\Big) + \PP\Big(Y'_h=1, B''_h=1, Y'_h\neq B'_h\Big). \notag 
\end{align}

We start by estimating the first term on the right hand side. Using H\"older's inequality, \eqref{term_1_1} and \eqref{term_1_3} we obtain for some $\gamma <1$ and $c<\infty$ depending on $h$ that
\begin{align*} 
\PP\Big(B'_h=1, Z''_h=1,E\Big) &\leq \PP\Big(B'_h=1, E\Big)^{\nicefrac 1 2} \cdot \PP\Big( Z''_h=1, E\Big)^{\nicefrac 1 2}\leq c \cdot \frac 1 n \cdot i \cdot \gamma^{i'+i''}.
\end{align*}
Since $i=i'+i''$ we can further enlarge $\gamma$ to incorporate $i$ and therefore obtain
\be
\PP\Big(B'_h=1, Z''_h=1,E\Big) \leq  c \cdot \frac 1 n \cdot \gamma^{i}.
\ee

The same argument works for evaluating the third, fifth and seventh term on the right hand side of \eqref{cov_1_SA_second_ineq} and gives the same bound.

In order to evaluate the second term on the right-hand side of \eqref{cov_1_SA_second_ineq} we follow the arguments from the proof of Lemma \ref{lemma_probs} used to show \eqref{term_1_2} in which we incorporate the event $\{Y''_h=1\}=\{B''_{\theta_h \wedge \theta_\tau}=1\}$ which is $\mathcal F_{\theta_h\wedge \theta_\tau}$-measurable. As in \eqref{intermediate} we obtain 
\begin{align*}
\PP\big( B'_h= 1, Y'_h \neq 1, Y''_h=1 \big) &  \leq \EE\Big((h-\theta_h\wedge \theta_\tau)\cdot B'_{\theta_h \wedge \theta_\tau}\cdot \Big( \frac {h}{h+2} \Big)^{B'_{\theta_h \wedge \theta_\tau}-2}; Y''_h=1 \Big)\end{align*}
which by H\"older's inequality, \eqref{prob_Y'} and using  further the arguments that lead to \eqref{term_1_2} gives
\be\label{second_term}
\PP\Big(B'_h=1,Y''_h=1, Y'_h\neq B'_h\Big) \leq c \cdot \frac 1 n \cdot \gamma^{i}
\ee
for some $\gamma <1$ and $c<\infty$ depending on $h$. We used the fact that $i'+i''=i<n$.

In order to bound the fourth term on the right-hand side of \eqref{cov_1_SA_second_ineq} we again a similar argument as in Lemma \ref{lemma_probs}, namely the one used to get \eqref{term_1_4}, where we incorporate the event $\{Y''_h=1\}=\{B''_{\theta_h \wedge \theta_\tau}=1\}$ and obtain the desired bound. 

For estimating the sixth term on the right-hand side of \eqref{cov_1_SA_second_ineq} observe that by H\"older's inequality
\begin{align*}
\PP\Big(B'_h=1, B''_h=1, Y'_h\neq B'_h\Big)  & \leq \PP\Big(B'_h=1, Y'_h\neq B'_h, Y''_h=1\Big) \\
& \qquad \qquad+\PP\Big(B'_h=1,Y'_h\neq B'_h, B''_h=1, Y''_h \neq B''_h\Big) 
\\
& \leq \PP\Big(B'_h=1, Y'_h\neq B'_h, Y''_h=1\Big) 
\\
& \qquad \qquad+\PP\Big(B'_h=1,Y'_h\neq B'_h\Big)^{\nicefrac 1 2} \cdot \PP\Big(B''_h=1, Y''_h \neq B''_h\Big)^{\nicefrac 1 2}
\end{align*}
which by means of \eqref{term_1_2} and \eqref{second_term}  can be again bounded from above by $c \cdot \frac 1 n \cdot \gamma^{i}$ for some $\gamma <1$ and $c<\infty$ depending on $h$. The same argument works also for the eighth term on the right-hand side of \eqref{cov_1_SA_second_ineq}. This completes the proof. 
\ep

\vspace{1cm}
To finish this section we deduce from the previous results that for $A' \cap A'' =\emptyset$ with $|A'|=i'$ and $|A''|=i''$ the following bound holds
\be\label{bound_covariances}
\Big| \mathbb {COV}\Big(\1_{\mathscr F_{A'}}, \1_{\mathscr F_{A''}}\Big)\Big| \leq c \cdot \frac 1 n \cdot \gamma^{i'+i''},
\ee
for some $\gamma <1$ and $c<\infty$ depending on $h$.

To see that this is a direct consequence of the previous three lemmas, observe that since the birth-death processes $B'$ and $B''$ are independent, as resulting from Lemma \ref{indep_B_j}, it holds that
\begin{align*}
\mathbb {COV}\Big(\1_{\mathscr F_{A'}}, \1_{\mathscr F_{A''}}\Big)& = \mathbb {COV}\Big(\1_{\{ Z'_h=1\}}, \1_{\{ Z''_h=1\}}\Big) 
\\  
& = \mathbb {COV}\Big(\1_{\{ B'_h=1\}}, \1_{\{ B''_h=1\}}\Big) \notag
\\  
&  \hphantom{12}+  \mathbb {COV}\Big(\1_{\{ B'_h=1\}},  \1_{\{ Z''_h=1\}}- \1_{\{ B''_h=1\}}\Big) + \mathbb {COV}\Big(\1_{\{ Z'_h=1\}}-\1_{\{ B'_h=1\}}, \1_{\{ Z''_h=1\}}\Big) \notag
\\  
& = \EE\Big(\1_{\{ B'_h=1\}} (\1_{\{ Z''_h=1\}}- \1_{\{ B''_h=1\}})\Big) - \PP(B'_h=1)\EE(\1_{\{ Z''_h=1\}}- \1_{\{ B''_h=1\}}) \notag
\\  
& \hphantom{12}+ \EE\Big( \1_{\{ Z''_h=1\}}(\1_{\{ Z'_h=1\}}-\1_{\{ B'_h=1\}})\Big) - \PP( Z''_h=1)\EE(\1_{\{ Z'_h=1\}}-\1_{\{ B'_h=1\}}) \notag
\\
& = \PP\Big(B'_h=1, Z''_h=1, Z''_h\neq B''_h\Big) - \PP\Big(B'_h=1, B''_h=1, Z''_h\neq B''_h\Big) \notag
\\  
& \hphantom{1234} - \PP(B'_h=1)\PP(Z''_h=1, Z''_h\neq B''_h) + \PP(B'_h=1)\PP(B''_h=1, Z''_h\neq B''_h)\notag
\\  
& \hphantom{12}+ \PP\Big(Z''_h=1, Z'_h=1, Z'_h\neq B'_h\Big) - \PP\Big(Z''_h=1, B'_h=1, Z'_h\neq B'_h\Big) 
\\  
& \hphantom{1234} - \PP( Z''_h=1)\PP( Z'_h=1, Z'_h\neq B'_h) + \PP( Z''_h=1)\PP( B'_h=1, Z'_h\neq B'_h).\notag
\end{align*}
Now \eqref{bound_covariances} follows directly by applying Lemma \ref{lemma_probs}, Lemma \ref{lemma_probs_2} and \eqref{BD_prob_i}.

\section{The initial and final levels of branches}\label{Lower and upper levels}

We start by making some notation. We fix a set $A \in \MP_n$, $A \neq \emptyset, \{1, \dots, n\}$, and view the coalescent tree of the population alive at time 0 from the leaves towards the root, that is in the coalescent time direction. Recall that $\Pi_k$ denotes the state of the coalescent after $n-k$ coalescing events. We denote the levels at which the branch with leaves in $A$ is formed and respectively ends by $K_A$ and $J_A$, thus 
\[K_A:=  \max\{ 2 \leq k \leq n: A \in \Pi_k\} \text{  \hspace{.3cm} and  \hspace{.3cm} } J_A:= \max\{ 1 \leq j < K_A: A \notin \Pi_j\}. \]
These are the {\em initial level} and the {\it final level} of the branch.  For a set $A$ of leaves which is not supported by some branch (which means that $A \notin \Pi_k$ for all $k$) it is convenient to set $K_A=J_A=\infty$. On the event $\{J_A<\infty\}=\{K_A<\infty\}$ the coalescent contains a subtree with root at level $J_A$, its last merger at level $K_A$ and its set of leaves equal to $A$. We denote this subtree by $\mathcal T_A$. The length of the branch supporting the leaves with labels in $A$ is
\be\label{def_SA}
L_A:=\sum_{j=J_A+1}^{K_A} X_j,
\ee
where the $X_j \sim Exp \Big(\binom j 2\Big)$ are the inter-coalescence times. For sets $A$ of leaves that are not supported by any branch, we have that $L_A=0$. Note that the whole tree structure, up to the coalescent times, is captured by the random variables $\{K_A, J_A\}_{A \in \MP_n\setminus \{\emptyset, \{1, \dots, n\}\}}$.

The following lemma gives bounds on the distribution weights of the initial and final levels of branches in the coalescent of the population alive at time 0. For the remainder of this section we consider the sets $A, A' \in \MP_n$ with $|A|=i$ and $|A'|=i'$ fixed, $1 \leq i, i'<n$.

\bl\label{lemma_J}
\begin{enumerate}[(i)]
\item For $1 \leq j < n$ and $1 \leq j' < n$ the following bounds hold:
\[\PP(J_A=j) \leq \frac{2 j }{ (n-1)}\cdot \frac 1 { \binom n i} \]
and for $A \neq A'$
\[ \PP(J_A=j, J_{A'}=j') \leq \begin{cases}\frac{4 j j' }{(n-1)(n-2)} \cdot \frac 1 {\binom n {i, i', n-(i+i')}} & \text{ if } A\cap A' =\emptyset, i+i'<n
						      \\[1.5 ex] \frac{ 4j(i-i') }{(n-1)(n-2)} \cdot \frac{1}{\binom n {i-i', i', n-i}} & \text{ if } A' \subset A \text{ and } j'> j.
			   	 \end{cases} 
\]
Otherwise $\PP(J_A=j, J_{A'}=j')=0$. 

Moreover, if $i+i'<n$ and $A\cap A' =\emptyset$ then for a finite constant $c$ depending on $i$ and $i'$
\[\PP(J_A=j, J_{A'}=j') \leq \PP(J=j)\PP(J'=j') \Big(1 +\frac c n\Big).\]

\item Let $i,i'\geq 2$. Then for $1 \leq j < k <n$ and $1 \leq j' <k' < n$ the following holds:
\[ \PP(K_A=k)\leq \frac{ik^2}{ (n-1)(n-2)}\cdot \frac 1 { \binom n i}, \quad \quad \PP(K_A=k, J_A=j)\leq \frac {2ij}{(n-1)(n-2)}\cdot \frac 1 { \binom n i}.\]

Moreover, if $i+i'< n$ and $A\cap A' =\emptyset$ then for a finite constant $c$ depending on $i$ and $i'$
\[\PP(K_A=k, J_A=j, K_{A'}=k', J_{A'}=j') \leq \PP(K_A=k, J_A=j)\PP(K_{A'}=k', J_{A'}=j') \Big(1 +\frac c n\Big).\]
\end{enumerate}
\el

\br\label{remark_ext}
Note that the bounds of the lemma are not valid in two exceptional cases: if $i=1$ then $K_A=n$ and if $i+i'=n$ and $A\cap A' =\emptyset$ then $J_A=1$ implies $J_{A'}=1$ such that
\begin{align*}
\PP(J_A=1, J_{A'}=1)& =  \PP(J_A=1)\leq  \frac{ 2}{(n-1) \cdot \binom{n} {i}}.
\end{align*}
\er

\bp
For this proof let $J:=J_A$, $J':=J_{A'}$, $K:=K_A$ and $K':=K_{A'}$.

(i) Note first that if $j>n-i$ then $\PP(J=j)=0$ and the claim holds trivially. Thus let $j \leq n-i$ and denote by $E_1$ the event that the subtree $\mathcal T_A$ has its $i-1$ mergers at the fixed levels $n> k_{1}  > \cdots > k_{i-1}> j$. To generate this event it is required that from level $n-1$ to level $j+1$ the branches from $\mathcal T_A$ merge at the levels $k_1, \dots, k_{i-1}$ and that at the remaining levels branches from outside of $\mathcal T_A$ merge. In addition, at level $j$ the remaining branch from $\mathcal T_A$ has to coalesce with one of the other $j$ branches extant in the coalescent tree. For the $i-1$ mergers within $\mathcal T_A$ there are $\binom{i}  {2} \cdot \binom{i-1}  {2}\cdots \binom {2} {2}$ available ways in which the branches can merge with one another, whereas for the other $n-i-j$ mergers there are $\binom{n-i}  {2} \cdot \binom{n-i-1} {2}\cdots \binom {j+1} {2}$ ways. Thus
\begin{align}\label{prob_E_1}
\PP(E_1)&= \frac {\binom{i}  {2} \cdot \binom{i-1}  {2}\cdots \binom {2} {2} \cdot \binom{n-i}  {2} \cdot \binom{n-i-1} {2}\cdots \binom {j+1} {2} \cdot j}{\binom{n}  {2} \cdot \binom{n-1} {2}\cdots \binom {j+1} {2}} \notag
\\
& = \frac{\binom{i}  {2}\cdots \binom {2} {2} \cdot j}{\binom{n}  {2} \cdot \binom{n-1} {2}\cdots \binom {n-i+1} {2}}
\\
& = \frac {2j}{\binom n i}  \cdot  \frac {(i-1)!}{(n-1)\cdots (n-i)}. \notag
\end{align}
Note that this probability does not depend on the choice of the levels $k_1, \dots, k_{i-1}$. There are $\binom {n-j-1}{i-1}$ possible ways to choose these levels and therefore
\begin{align}\label{distrib_J_1}
\PP(J=j)= \binom {n-j-1}{i-1} \cdot \PP(E_1) &= \frac {2j}{\binom n i} \cdot \frac {(n-j-1)!} {(n-1)\cdots (n-i)(n-j-i)! }
\leq\frac{2 j }{ (n-1)}\cdot \frac 1 { \binom n i}. 
\end{align}

Let us now look at the joint distribution of $J$ and $J'$ and consider first the case where $A \cap A' = \emptyset$ and $i+i' < n$. Let us assume first that $j'<j$. If $j>n-i$ or $j'>n-i-i'$ then $\PP(J=j, J'=j')=0$. Otherwise, the $i'-1$ levels at which mergers in $\mathcal T_{A'}$ take place are all different from the $i-1$ levels at which the mergers in $\mathcal T_A$ occur and they are also different from $j$. At these $i-1$ (and respectively $i'-1$) levels pairs of branches from the subtree $\mathcal T_A$ (respectively from $\mathcal T_{A'}$) are chosen to merge, whereas at the other levels (from level $n-1$ to level $j'+1$, excepting $j$) branches with leaves in $\{1, \dots, n\} \setminus (A \cup A')$ merge. At level $j$ the branch with leaves in $A$ ends through a merger with one of the $j-q$ existing branches with leaves in $\{1, \dots, n\} \setminus (A \cup A')$. At level $j'$ the branch with leaves in $A'$ ends by merging with one of the $j'$ branches with leaves in $\{1, \dots, n\} \setminus (A \cup A')$. Similar as above denote by $E_2$ the event that the subtrees $\mathcal T_A$ and $\mathcal T_{A'}$ have their mergers at the fixed levels $n> k_{1}  > \cdots > k_{i-1}> j$ and $n> k'_{1}  > \cdots > k'_{i'-1}> j$ respectively. We obtain
\begin{align}\label{prob_E_2}
\PP(E_2)&= \frac {\binom{i}  {2} \cdots \binom {2} {2} \cdot \binom{i'} {2}\cdots \binom {2} {2}  \cdot \binom{n-(i+i')} {2} \cdots \binom {j'+1} {2} \cdot (j-q) \cdot j'}{\binom{n}  {2} \cdot \binom{n-1} {2}\cdots \binom {j'+1} {2}} \notag
\\
& \leq \frac{\binom{i}  {2}\cdots \binom {2} {2} \cdot \binom{i'}  {2}\cdots \binom {2} {2}\cdot j j'}{\binom{n}  {2} \cdot \binom{n-1} {2}\cdots \binom {n-(i+i')+1} {2}}
\\
& \leq \frac {4jj'}{\binom n {i, i', n-(i+i')}}  \cdot  \frac {(i-1)!(i'-1)!}{(n-1)\cdots (n-(i+i'))}. \notag
\end{align}
There are $\binom{n-j-1}{i-1} \binom {n-j'-i-1}{i'-1}$ possible ways to first choose the levels $k_1, \dots, k_{i-1}$ and then $k'_{1}, \dots, k'_{i'-1}$ such that $\{k'_1, \dots, k'_{i'-1} \}  \cap \{k_1, \dots, k_{i-1}, j\}=\emptyset$. Thus, since $1\leq j' < j$,
\begin{align}\label{joint_distrib_J}
\PP(J=j, J'=j')&= \binom{n-j-1}{i-1} \binom {n-j'-i-1}{i'-1} \cdot \PP(E_2) 
\\
& \leq \frac{4 j j'} {\binom n {i, i', n-(i+i')}} \notag
\\
& \quad\quad \cdot \frac{(n-j-1)!}{(n-1)\cdots(n-i)(n-j-i)!} \cdot \frac{(n-j'-i-1)!}{(n-i-1)\cdots(n-i-i')(n-j'-i-i')!}\notag
\\
& \leq \frac{4 j j' }{(n-1)(n-2)} \cdot \frac 1 {\binom n {i, i', n-(i+i')}}. \notag
\end{align}

Observe also that by putting together \eqref{joint_distrib_J}, \eqref{prob_E_2}, \eqref{prob_E_1} and \eqref{distrib_J_1} we get
\begin{align*}
\PP(J=j, J'=j')&\leq \binom{n-j-1}{i-1}  \frac{ \binom{i} {2} \cdots \binom{2} {2} }{\binom{n} {2} \cdots \binom {n-i+1} {2}}\cdot j\cdot \binom {n-j'-1}{i'-1} \frac{\binom{i'} {2} \cdots \binom{2} {2} }{\binom{n-i} {2} \cdots \binom {n-(i+i')+1} {2}} \cdot j'
\\
&= \PP(J=j)\PP(J'=j')\cdot \frac{\binom{n} {2} \cdots \binom {n-i'+1} {2} }{\binom{n-i} {2} \cdots \binom {n-i-i'+1} {2}}
\\
& \leq \PP(J=j)\PP(J'=j') \Big(1+\frac c n\Big)
\end{align*}
for a finite constant $c$ depending on $i$ and $i'$.

In the case $j=j'$, $i+i'<n$ we have that $j=j'\geq 2$. If $j > n-i-i'+1$ then $\PP(J=j, J'=j)=0$. Otherwise the two branches we consider end at level $j$ by coalescing with one another and therefore we obtain with a similar argument as above that 
\begin {align}\label{joint_distrib_J_2}
\PP(J=j, J'=j) &=  \binom {n-j-1} {i-1}\binom{n-j-i}{i'-1}   \frac {\binom{i}  {2} \cdots \binom {2} {2} \cdot \binom{i'} {2}\cdots \binom {2} {2}  \cdot \binom{n-(i+i')} {2} \cdots \binom {j} {2} \cdot 1}{\binom{n}  {2} \cdot \binom{n-1} {2}\cdots \binom {j+1} {2}}  \notag
\\  
& = \binom {n-j-1} {i-1}\binom{n-j-i}{i'-1}   \frac{ \binom{i} {2} \cdots \binom{2} {2}\cdot  \binom{i'} {2} \cdots \binom{2} {2} }{\binom{n} {2} \cdots \binom {n-(i+i')+1} {2}}  \cdot \binom j 2 
\\  
&= \frac {4 \binom j 2} {\binom n {i, i', n-(i+i')}}\cdot \frac 1 {(n-1)(n-2)} \cdot \frac{(n-j-1)\cdots (n-j-i-i'+2)}{(n-3)\cdots(n-i-i')}. \notag
\end{align}
Using $j \geq 2, i+i' <n$ we get that
\begin{align*}
\PP(J=j, J'=j)\leq \frac{2 j^2}{(n-1)(n-2)} \cdot \frac 1 {\binom n {i, i', n-(i+i')}}. \notag
\end{align*}

Also in this case by putting together \eqref{joint_distrib_J_2}, \eqref{prob_E_1} and \eqref{distrib_J_1} we get
\begin{align*}
\PP(J=j, J'=j)&\leq \binom{n-j-1}{i-1}\frac{\binom{i} {2} \cdots \binom{2} {2}}{\binom{n} {2} \cdots \binom {n-i+1} {2}}\cdot j \cdot \binom{n-j-1}{i'-1}  \frac{\binom{i'} {2} \cdots \binom{2} {2} }{\binom{n-i} {2} \cdots \binom {n-(i+i')+1} {2}} \cdot \frac{j-1} 2
\\
& \leq \PP(J=j)\PP(J'=j) \cdot \frac{\binom{n} {2} \cdots \binom {n-i'+1} {2} }{\binom{n-i} {2} \cdots \binom {n-i-i'+1} {2}} 
\\
& \leq \PP(J=j)\PP(J'=j) \Big(1+\frac c n\Big)
\end{align*}
for a finite constant $c$ depending on $i$ and $i'$.

Let now $A' \subset A$, $A \neq A'$. For $j'>n-i'$ or $j>n-i$ it holds that $\PP(J=j, J'=j)=0$. Otherwise $\mathcal T_{A'}$ is a subtree of $\mathcal T_A$ and $j'> j$. If $\mathcal T_A$ has its mergers at the levels $k_1, \dots, k_{i-1}$, then $j' \in \{k_1, \dots, k_{i-1}\}$ and at level $j'$ the remaining branch of $\mathcal T_{A'}$ merges with one of the $q$ existing branches with leaves in $A\setminus A'$, where $1 \leq q \leq i-i'$. By similar arguments as above we deduce that
\begin{align*}
\PP(J=j, J'=j') \leq \binom{n-j'-1}{i'-1}\binom{n-j-i'-1}{i-i'-1}\frac{ \binom{i-i'} {2} \cdots \binom{2} {2}\cdot  \binom{i'} {2} \cdots \binom{2} {2}  }{\binom{n} {2} \cdots \binom {n-i+1} {2}} \cdot  j \cdot (i-i').
\end{align*}
Replacing $j$ by 1 and $j'$ by 2 in the binomials we obtain
\begin{align*}
\PP(J=j, J'=j') & \leq 4j(i-i') \cdot\frac{1}{(n-1)(n-2)} \cdot \frac{1}{\binom n {i-i', i', n-i}}.
\end{align*}

\vspace{0.8cm}
(ii) Let now $i \geq 2$. If $k >n-i$ then $\PP(K=k, J=j) =\PP(K=k) =0$. Otherwise the probability that the subtree $\mathcal T_A$ has its last merger at level $k$ and its root at level $j$ can be quickly obtained by summing up the equality \eqref{prob_E_1} over all possible ways of choosing the merging levels $n> k_1> \dots> k_{i-2}>k$. Thus
\begin{align}\label{proof_K_J}
\PP(K=k, J=j) = \binom {n-k-1} {i-2} \frac{\binom{i}  {2} \cdots \binom {2} {2}}{\binom{n} {2} \cdots \binom{n-i+1} {2} } \cdot j \leq \frac {2ij}{(n-1)(n-2)}\cdot \frac 1 { \binom n i}. 
\end{align}
Summing now over all $1 \leq j<k$ we get that
\begin {align*}
\PP(K=k)& \leq \sum_{j=1}^k \frac {2ij}{(n-1)(n-2)}\cdot \frac 1 { \binom n i} \leq \frac {ik^2}{(n-1)(n-2)}\cdot \frac 1 { \binom n i}.
\end{align*}

As to the last claim of the lemma, $\PP(K=k, J=j, K'=k', J'=j')=0$ if $k=k'$. Without loss of generality assume $k'<k$. Note that the condition $i+i'<n$ implies that $J>1$ or $J'>1$. Considering the positions of $j$, $j'$ and $k'$ we have to distinguish between the cases $j<j'$, $j=j'$, $j'<j<k'$ and $k'<j$. In the first case, one needs to choose the levels $k_1, \dots, k_{i-2}$ and $k'_1, \dots, k'_{i'-2}$ such that 
\[ \{k'_1, \dots, k'_{i'-2} \}  \cap \{k_1, \dots, k_{i-2}, k\}=\emptyset. \]
Since $j<j'$ the branch supporting $\mathcal T_{A'}$ (which ends at level $j'$) cannot merge with the branch supporting $\mathcal T_A$ and must therefore merge with one of the $j'-1$ branches left. Thus
\begin{align}\label{case_1_intermed}
\PP(K=k, J=j, K'=k', J'=j')&=  \binom {n-k-1} {i-2}  \binom {n-k'-i} {i'-2} \frac{\binom{i}  {2} \cdots \binom {2} {2}\cdot \binom{i'}  {2} \cdots \binom {2} {2}}{\binom{n} {2} \cdots \binom{n-(i+i')+1} {2} } \cdot j\cdot (j'-1)
\end{align}
which by \eqref{proof_K_J} leads to
\begin{align*}
\PP(K=k, J=j, K'=k', J'=j')&\leq \PP(K=k, J=j)\PP(K'=k', J'=j')\cdot \frac{\binom{n} {2} \cdots \binom {n-i'+1} {2} }{\binom{n-i} {2} \cdots \binom {n-(i+i')+1} {2}}
\\ 
&\leq \PP(K=k, J=j)\PP(K'=k', J'=j')  \Big(1+\frac c n\Big)
\end{align*}
for a finite constant $c$ depending on $i$ and $i'$.

In the case $j'<j<k'$ the same argument holds, with the observation that the product $j\cdot (j'-1)$ on the right-hand side of \eqref{case_1_intermed} has to be replaced by $(j-1) \cdot j'$, whereas in the case $j=j'$ it is $\binom j 2$ that appears instead of the product $j\cdot (j'-1)$. In the last case, $k'<j$, observe that the intermediate levels $k_1, \dots, k_{i-2}$ and $k'_1, \dots, k'_{i'-2}$ have to be chosen such that
\[ \{k'_1, \dots, k'_{i'-2} \}  \cap \{k_1, \dots, k_{i-2}, k,j\}=\emptyset \]
and thus
\begin{align*}
\PP(K=k, J=j, K'=k', J'=j')&=  \binom {n-k-1} {i-2}  \binom {n-k'-i-1} {i'-2} \frac{\binom{i}  {2} \cdots \binom {2} {2}\cdot \binom{i'}  {2} \cdots \binom {2} {2}}{\binom{n} {2} \cdots \binom{n-(i+i')+1} {2} } \cdot (j-q)\cdot j'
\end{align*}
where $q$ denotes as above the number of branches supporting leaves with labels in $A'$ extant at level $j$ (and with which the branch supporting the leaves with labels in $A$ can not merge). By \eqref{proof_K_J} we obtain the desired upper bound. This finishes the proof.
\ep

The following lemma employs Lemma \ref{lemma_J} extensively and will be of use in the proof of Theorem \ref{mainresult_fdd}.

\bl\label{lemma_claims}
For fixed $i\geq 1$ as $n \to \infty$
\begin{equation}\label{claim_1}
\frac {n}{\log n} \sum_{A\in \MP_n \atop |A|=i} \Big(\frac 1 {J_A} - \EE\Big(\frac 1 {J_A}\Big)\Big)^2 \stackrel {\PP} \longrightarrow 1
\end{equation}
with $\frac 1 {J_A}=0$ if $J_A=\infty$. Also
\begin{equation}\label{claim_2}
\max_{A\in \MP_n \atop |A|=i} \sqrt{\frac {n}{\log n} }\Big|\frac 1 {J_A} - \EE\Big(\frac 1 {J_A}\Big)\Big|\stackrel {\PP} \longrightarrow 0.
\end{equation}
\el

\bp
Remember that by definition $J_A=\infty$ if there exists no branch in the coalescent supporting the leaves with labels in the set $A$ and otherwise $J_A<n$.  It thus holds that
\begin{align*}
\frac {n}{\log n} \sum_{A\in \MP_n \atop |A|=i}\frac 1 {J_A^2} & = \frac {n}{\log n} \sum_{A\in \MP_n \atop |A|=i}\frac 1 {J_A^2}\cdot \1_{\{J_A<n\}}. \notag
\end{align*}
In what follows we show that 
\begin{equation}\label{conv_J}
\frac {n}{\log n} \sum_{ |A|=i} \frac 1 {J_A^2}\cdot \1_{\{J_A<n\}} \stackrel {\PP} \longrightarrow 1.
\end{equation}

Let us denote by $\mathcal A=\mathcal A(n,i ,a)$ the event that at least one branch of order $i$ ends between level 2 and level $a \sqrt n$ in the coalescent, where $a$ is a positive constant:
\[\mathcal A=\mathcal A(n, i,a):=\{ J_A < a \sqrt n \quad \text{ for some } A\in \MP_n, |A|=i \}.\]
Note by Lemma \ref{lemma_J} that
\begin{equation}\label{A_c}
\PP(\mathcal A) \leq \sum_{ |A|=i} \PP(J_A < a \sqrt n) \leq  \sum_{ |A|=i} \sum_{1 \leq j<a\sqrt n} \frac{2 j }{ (n-1)}\cdot \frac 1 {\binom n i} \leq 2a^2.
\end{equation}
Since $ \sum_{ |A|=i} \frac 1 {J_A^2}\cdot \1_{\{J_A<n\}} \neq  \sum_{ |A|=i} \frac 1 {J_A^2}\cdot \1_{\{a\sqrt n\leq J_A <n\}}$ only on the event $\mathcal A$, which for small $a$ has arbitrary small probability, in order to prove \eqref{conv_J} it suffices to show that for all $a>0$
\[ \frac {n}{\log n} \sum_{ |A|=i} \frac 1 {J_A^2}\cdot \1_{\{a\sqrt n\leq J_A <n\}}  \stackrel {\PP} \longrightarrow 1 \quad \text{as $n \to \infty$.} \] 

Using \eqref{distrib_J_1} we obtain that for all $a >0$
\begin{align}\label{expect_J}
\EE\Big(\frac {n}{\log n} \sum_{ |A|=i} \frac 1 {J_A^2}\cdot \1_{\{a\sqrt n\leq J_A <n\}}\Big)& =\frac {n}{\log n} \sum_{ |A|=i} \sum_{a\sqrt n \leq j<n}  \frac 1 {j^2} \cdot \frac{2j}{\binom n i} \cdot\frac {(n-j-1)\cdots (n-j-i+1)}{(n-1) \cdots (n-i)} \notag
\\
& =\frac {n}{(n-1)\log n} \sum_{a\sqrt n \leq j<n} \frac 2 {j} \cdot \frac {(n-2-(j-1))\cdots (n-i-(j-1))}{(n-2) \cdots (n-i)} \notag
\\
& \to 1 
\end{align}
as $n \to \infty$. 

For the variance we obtain
\begin{align}\label{var_J}
\VV\Big(\frac {n}{\log n} \sum_{ |A|=i} \frac 1 {J_A^2}\cdot \1_{\{a\sqrt n\leq J_A <n\}}\Big)&\leq \frac {n^2}{\log^2 n}\sum_{ A} \EE\Big(\frac {1} {J_A^4}\cdot \1_{\{a\sqrt n\leq J_A <n\}} \Big)
\\
& \quad + \frac {n^2}{\log^2 n}\sum_{A \cap A' = \emptyset} \mathbb{COV}\Big(\frac 1 {J_A^2}\cdot \1_{\{a\sqrt n\leq J_A <n\}}, \frac 1 {J_{A'}^2}\cdot \1_{\{a\sqrt n\leq J_{A'} <n \}}\Big). \notag
\end{align}
Note that for $A \cap A' \notin \{A, A', \emptyset\}$, $|A|=|A'|=i$ either $J_A$ or $J_{A'}=\infty$ and thus the corresponding covariances are less than or equal to 0.

By Lemma \ref{lemma_J} it holds for $A \cap A' = \emptyset$ that
\begin{align*}
\mathbb{COV}\Big(\frac 1 {J_A^2}\cdot \1_{\{a\sqrt n\leq J_A <n\}}, &\frac 1 {J_{A'}^2}\cdot \1_{\{a\sqrt n\leq J_{A'} <n\}}\Big)
\\
& = \sum_{a\sqrt n \leq j,j'<n} \frac {1}{j^2 j'^2} \cdot \Big(  \PP(J_{A}=j, J_{A'}=j') - \PP(J_{A}=j) \PP(J_{A'}=j')  \Big)
\\
& \ll \sum_{a\sqrt n \leq j,j'<n} \frac {1}{j^2 j'^2} \cdot \PP(J_{A}=j) \PP(J_{A'}=j')  \cdot \frac 1 n.
\\
& \ll \sum_{a\sqrt n \leq j,j'<n}\frac {1}{j j'} \cdot \frac 1 {n^2} \cdot \frac 1{\binom n i ^2} \cdot \frac{1}{n}  
\\
& \ll \frac{(\log n)^2}{n^{3}}\cdot \frac 1{\binom n i ^2},
\end{align*}
for a constant $c$ depending on $i$ and $i'$. Plugging this in \eqref{var_J} and using again Lemma \ref{lemma_J} we obtain that
\begin{align*}
\VV\Big(&\frac {n}{\log n} \sum_{ |A|=i} \frac 1 {J_A^2}\cdot \1_{\{a\sqrt n\leq J_A <n\}}\Big)
\\
&\ll \frac {n^2}{(\log n)^2} \sum_{ A} \sum_{a\sqrt n \leq j<n} \frac {1} {j^4} \cdot \frac{2j}{n-1} \cdot \frac 1 {\binom n i}+ \frac {n^2}{(\log n)^2}\sum_{A \cap A' = \emptyset} \frac{(\log n)^2}{n^{3}}\cdot \frac 1{\binom n i ^2}
\\
& \ll \frac {n^2}{(\log n)^2} \Big( \frac 1 {a^2 n^2} + \frac{(\log n)^2}{n^{3}}\Big)
\\
& \ll \frac 1 {(\log n)^2} + \frac{1}{n}
\end{align*}
and thus
\begin{align*}
\VV\Big(\frac {n}{\log n} \sum_{|A|=i} \frac 1 {J_A^2}\cdot \1_{\{a\sqrt n\leq J_A <n\}}\Big)\to 0 \quad\quad \text{as } n\to \infty.
\end{align*}
This together with \eqref{expect_J} and \eqref{A_c} proves \eqref{conv_J}.

From Lemma \ref{lemma_J}
\be\label{expect_1/J}
\EE\Big(\frac 1 {J_A}\Big) \leq \sum_{j=1}^{n-1} \frac 1 j \cdot \frac {2j}{(n-1)} \cdot \frac 1 {\binom n i} = \frac 2 {\binom n i}
\ee
and thus for $i \geq 1$
\begin{align*}
\frac {n}{\log n} \sum_{ |A|=i} \EE\Big(\frac 1 {J_A}\Big)^2\ll \frac {n}{\log n} \sum_{A} \frac 1 {{\binom n i}^2} = \frac {n}{\log n} \cdot \frac 1 {\binom n i} \to 0
\end{align*}
as $n \to \infty$. Because of \eqref{conv_J} it follows by the Cauchy-Schwarz inequality that 
\[\frac {n}{\log n} \sum_{ |A|=i} 2\Big(\frac 1 {J_A}\Big)\cdot \EE\Big(\frac 1 {J_A}\Big) \leq \frac {2n}{\log n} \sqrt{ \sum_{A} \frac 1 {J_A^2}} \cdot \sqrt{\sum_{A} \EE\Big(\frac 1 {J_A}\Big)^2 }\stackrel {\PP} \longrightarrow 0 \]
as $n\to \infty$ and thus \eqref{claim_1} holds.

Finally observe that by \eqref{expect_1/J} 
\[\Big|\frac 1 {J_A}- \EE\Big(\frac 1 {J_A}\Big)\Big| \leq \frac 1 {J_A}+ \frac 2 n \leq \frac 1 {a \sqrt n} + \frac 2 n\]
on the event $\mathcal A^c$. Thus the claim \eqref{claim_2} is a direct consequence of \eqref{A_c}.
\ep

\section{Contribution of big families}\label{Proof of the Proposition}

In this section we show that the branches of big orders from the coalescent at time 0 make only a negligible contribution to the external length at a later time $h$ when the number of individuals in the population gets large. 

Recall the definition of $R^{n,r}_h$ given in \eqref{def_R} and let
\be\label{def_a_nh}
a_{n, h}:=\Big(\log\Big(\frac {h+2} h\Big)\Big)^{-1} \log n.
\ee
Then we can write 
\[
R^{n,r}_h= \widetilde R^{n,r}_h+\bm\bar R^{n,r}_h,
\]
where
\[\widetilde R^{n,r}_h:= \sum_{A \in \MP_n \atop r<|A|\leq a_{n,h} } \1_{\mathscr F_A } \cdot L_A, \qquad \qquad \bm\bar R^{n,r}_h:= \sum_{A \in \MP_n \atop a_{n,h}<|A|< n } \1_{\mathscr F_A} \cdot L_A .\]

\bprop\label{proposition_bigger_r}
For each $h \geq 0$ there exists a sequence $\{\varepsilon(r)\}_{r \in \NN}$ with $\varepsilon(r) = o(1)$ when $r \to \infty$, such that
\[\mathbb{V}\Big(\widetilde R^{n,r}_h\Big) \leq \varepsilon(r) \, \frac{\log n}{n} \] 
holds. Moreover $\EE\big(\bm\bar R^{n,r}_h\big)=O\big(\frac 1 n\big)$.
\eprop

\bp
First we prove the second claim. Since the collection of random variables $\{\1_{\mathscr F_A}\}_{A \in \MP_n}$ is independent of the collection of random variables $\{L_A\}_{A \in \MP_n}$ due to the Poissonian structure embedded in the Moran model, it holds that
\[
\EE \Big( \sum_{ |A|>  a_{n, h} } \1_{\mathscr F_A} \cdot L_A \Big) = \sum_{i >  a_{n, h} } \sum_{ |A|=i}\EE( \1_{\mathscr F_A}) \cdot \EE(L_A).
\]
Note that with the notation introduced in Section \ref{A birth and death process} it holds that 
\[\mathscr F_A=\{Z^A_h=1\}\] 
and thus by \eqref{expectation_1_SA} for $|A|=i$
\[p_{i,h}^n:=\PP(\mathscr F_A) \leq 2i \cdot \Big(\frac h {h+2}\Big)^{i-1}.\] 
By \eqref{expectation_L} it holds for $|A|=i$
\be\label{expectation_L_A}
\EE(L_A)=  \frac 1{\binom n i} \cdot\EE(\ML^{n,i}) = \frac 1{\binom n i}\cdot \frac 2 i.
\ee
We obtain that
\begin{align*}
\EE\big(\bm\bar R^{n,r}_h\big)=\EE \Big(\sum_{ |A|>  a_{n, h} } \1_{\mathscr F_A} \cdot L_A  \Big)
& \leq   \sum_{i >  a_{n, h} } \binom n i \cdot 2i \cdot \Big(\frac h {h+2}\Big)^{i-1} \cdot \frac 2 i \cdot \frac 1 {\binom n i}
\\ 
& \ll \Big(\frac h {h+2}\Big)^{ \left(\log\left(\frac {h+2} h\right)\right)^{-1} \log n} 
\\ 
& = \frac{1}{n}.
\end{align*}
Thus $\EE\big(\bm\bar R^{n,r}_h\big)=O\big(\frac 1 n\big)$.

\vspace{0.5cm}

We now proceed to proving the first claim of the Proposition. For the variance of $\widetilde R^{n,r}_h$ it holds that
\begin{align*}
\VV\Big(\widetilde R^{n,r}_h\Big) & = \VV\Big(\sum_{ r<|A|\leq a_{n,h}} \1_{\mathscr F_A } \cdot L_A\Big) 
\\ 
& = \VV\Big(\sum_{ r<|A|\leq a_{n,h}} (\1_{\mathscr F_A } - \EE(\1_{\mathscr F_A })) \cdot L_A+ \sum_{ r<|A|\leq a_{n,h}} \EE(\1_{\mathscr F_A }) \cdot L_A\Big) 
\\ 
& = \VV\Big(\sum_{ r<|A|\leq a_{n,h}} (\1_{\mathscr F_A } - \EE(\1_{\mathscr F_A })) \cdot L_A \Big)+ \VV\Big(\sum_{ r<|A|\leq a_{n,h}} \EE(\1_{\mathscr F_A }) \cdot L_A\Big) 
\\ 
& \hphantom{12}+ \mathbb{COV}\Big(\sum_{ r<|A|\leq a_{n,h}} (\1_{\mathscr F_A } - \EE(\1_{\mathscr F_A })) \cdot L_A, \sum_{ r<|A|\leq a_{n,h}} \EE(\1_{\mathscr F_A }) \cdot L_A\Big).
\end{align*}
Taking into account the independence of the collections of random variables $\{\1_{\mathscr F_A }\}_{A \in \MP_n}$ and $\{L_A\}_{A \in \MP_n}$, the last term on the right-hand side above is equal to zero. Therefore, using the fact that $\EE(\1_{\mathscr F_A })$ depends on the set $A$ only through its cardinality $i$ and again from independence, we have that
\begin{align}\label{variance_M}
\VV\Big(\widetilde R^{n,r}_h\Big) & = \EE\Big(\Big(\sum_{ r<|A|\leq a_{n,h}} (\1_{\mathscr F_A } - \EE(\1_{\mathscr F_A })) \cdot L_A\Big)^2 \Big)+ \VV\Big(\sum_{ r<|A|\leq a_{n,h}} \EE(\1_{\mathscr F_A }) \cdot L_A\Big) 
\\ 
& = \sum_{ r<|A|, |A'|\leq a_{n,h}}\EE\Big((\1_{\mathscr F_A } - \EE(\1_{\mathscr F_A }))(\1_{\mathscr F_{A' }} - \EE(\1_{\mathscr F_{A'} }))\Big) \cdot \EE( L_A L_{A'}) + \VV\Big(\sum_{ r<i \leq a_{n,h}} p_{i,h}^n \cdot \ML^{n,i}\Big). \notag
\end{align}
We start by evaluating the first term on the right-hand side. To this aim note that if $A$ and $A'$ are such that $A \cap A' \notin \{\emptyset, A, A'\}$, then there cannot be two branches in the coalescent tree, one supporting the leaves in $A$ and the other one supporting the leaves in $A'$ and therefore in this case the product $L_A L_{A'}$ is equal to zero.  Thus, suppressing in the notation the restrictions on the cardinalities of the sets, it holds that
\begin{align}\label{M_first_term_1}
\sum_{A, A' \in \MP_n }\EE\Big((\1_{\mathscr F_A } - &\EE(\1_{\mathscr F_A }))(\1_{\mathscr F_{A' }} - \EE(\1_{\mathscr F_{A'} }))\Big) \cdot \EE( L_A L_{A'})   \notag
\\ 
& \leq 2 \sum_{A' \varsubsetneqq A }\EE\Big(\1_{\mathscr F_A }\1_{\mathscr F_{A' }} \Big) \cdot \EE( L_A L_{A'})
\\
&\qquad \qquad + \sum_{A \cap A' = \emptyset}\mathbb{COV}\Big(\1_{\mathscr F_A },\1_{\mathscr F_{A' }} \Big) \cdot \EE( L_A L_{A'}) + \sum_{A}\VV\Big(\1_{\mathscr F_A }\Big) \cdot \EE( L_A^2). \notag
\end{align}

Using the definition (\ref{def_SA}) of $L_A$ we obtain that
\begin{align*}
\EE(L_A L_{A'}) &= \EE\Big(\EE(L_A L_{A'} \mid J_A, J_{A'}, K_A, K_{A'})\Big) 
\\ 
& = \EE\Big(\EE\Big(\sum_{j=J_A+1}^{K_A} X_j \cdot \sum_{l=J_{A'}+1}^{K_{A'}} X_l \mid J_A, J_{A'}, K_A, K_{A'}\Big)\Big)
\\ 
& = \EE\Big(\sum_{j=J_A+1}^{K_A} \sum_{l=J_{A'}+1}^{K_{A'}} \EE(X_j X_l)\Big).
\end{align*}
The exponential inter-coalescence times $X_j$ and $X_l$ are independent if $j \neq l$ and in this case $ \EE(X_j X_l) = \frac 2 {j(j-1)} \frac 2 {l(l-1)}$, whereas if $j=l$ it holds that $\EE(X_j X_l)=\EE(X_j^2)= 2 \, \Big( \frac 2 {j(j-1)}\Big)^2 =  2 \,\frac 2 {j(j-1)} \frac 2 {l(l-1)} $. We thus obtain that
\begin{align}\label{expectation_SA_SA'}
\EE(L_A L_{A'}) &\leq 2 \, \EE\Big(\sum_{j=J_A+1}^{K_A} \sum_{l=J_{A'}+1}^{K_{A'}} \frac 2 {j(j-1)} \frac 2 {l(l-1)}  \Big) \notag
\\ 
& = 2 \, \EE\Big(\sum_{j=J_A+1}^{K_A}  \frac 2 {j(j-1)} \cdot \sum_{l=J_{A'}+1}^{K_{A'}} \frac 2 {l(l-1)}  \Big) \notag
\\ 
& = 2 \,  \EE\Big(\Big( \frac 2 {J_A} -\frac 2 {K_A} \Big)\Big( \frac 2{J_{A'}}- \frac 2{K_{A'}}\Big) \Big) \notag
\\ 
&\leq 8 \, \EE\Big( \frac 1 {J_AJ_{A'}} \Big). 
\end{align}

In order to bound the first term on the right-hand side of (\ref{M_first_term_1}) observe first that by (\ref{expectation_1_SA})
\[\EE\Big(\1_{\mathscr F_A }\1_{\mathscr F_{A' }} \Big) \leq \EE\Big(\1_{\mathscr F_A }\Big) \leq 2i \Big(\frac h {h+2}\Big)^{i-1}.\]
Now, using (\ref{expectation_SA_SA'}) together with the fact that $A' \varsubsetneqq A$ implies that $J_{A'} >J_{A}$ and Lemma \ref{lemma_J} we obtain

\begin{align}\label{M_first_term_2}
&\sum_{A' \varsubsetneqq A }\EE\Big(\1_{\mathscr F_A }\1_{\mathscr F_{A' }} \Big) \cdot \EE( L_A L_{A'}) \notag
\\ 
& \leq 16 \, \sum_{i'<i} \binom n {i-i', i', n-i} \cdot i \Big(\frac h {h+2}\Big)^{i-1} \cdot \sum_{1 \leq j < j'<n} \frac 1 {jj'} \,\PP(J_A=j, J_{A'}=j') \notag
\\ 
&\leq 16 \, \sum_{i'<i} \binom n {i-i', i', n-i} \cdot i \Big(\frac h {h+2}\Big)^{i-1} \cdot \sum_{1 \leq j < j'<n} \frac 1 {jj'} \, \frac{ 4j(i-i') }{(n-1)(n-2)} \cdot \frac{1}{\binom n {i-i', i', n-i}} \notag
\\ 
&\leq 64 \, \sum_{i'<i} i \Big(\frac h {h+2}\Big)^{i-1} \cdot \frac{ i-i' }{(n-1)(n-2)} \sum_{1 \leq j< j'<n} \frac 1 {j'} \notag
\\ 
&\ll \frac{ 1 }{n} \cdot \sum_{i'<i} i ( i-i' )\Big(\frac h {h+2}\Big)^{i-1} \notag
\\ 
&\ll \frac{ 1 }{n} \cdot \sum_{i>r} i^3 \Big(\frac h {h+2}\Big)^{i-1} \notag
\\ 
&\ll \frac{ 1 }{n}. 
\end{align}

For bounding the second term on the right-hand side of (\ref{M_first_term_1}) observe that for $A \cap A' =\emptyset$ by (\ref{expectation_SA_SA'}) and Lemma \ref{lemma_J} it holds that
\begin{align*}
\EE(L_A L_{A'}) & \leq 8 \, \sum_{j, j'\geq 1} \frac 1 {jj'} \, \PP(J_A=j, J_{A'}=j') 
\\ 
& \leq  32 \, \sum_{j, j'\geq 1} \frac 1 {jj'} \, \frac{ j j' }{(n-1)(n-2)} \cdot \frac 1 {\binom n {i, i', n-(i+i')}} \notag
\\ 
& \ll \frac 1 {\binom n {i, i', n-(i+i')}} \notag
\end{align*}
and therefore by \eqref{bound_covariances} we obtain that for some $\gamma <1$ depending on $h$ it holds
\begin{align}\label{covs}
 \sum_{A \cap A' = \emptyset}\mathbb{COV}\Big(\1_{\mathscr F_A },\1_{\mathscr F_{A' }} \Big) \cdot \EE(L_A L_{A'}) & \ll \sum_{A \cap A' = \emptyset} \frac 1 n \cdot \gamma^{i+i'} \cdot \frac 1 {\binom n {i, i', n-(i+i')}} \ll \frac 1 n \sum_{i,i'>r} \gamma^{i+i'} \ll \frac{ 1 }{n}. 
\end{align}

We now turn to the last term on the right-hand side of (\ref{M_first_term_1}), namely $\sum_{A}\VV\Big(\1_{\mathscr F_A }\Big) \cdot \EE( L_A^2)$. From (\ref{expectation_1_SA}) it holds that
\[
\VV\Big(\1_{\mathscr F_A } \Big) \leq \EE\Big(\1_{\mathscr F_A }\Big) \leq 2i \Big(\frac h {h+2}\Big)^{i-1}
\]
and by (\ref{expectation_SA_SA'}) and Lemma \ref{lemma_J} 
\begin{align*}
\EE(L_A^2) \leq 8 \, \sum_{1 \leq j<n} \frac 1 {j^2} \, \PP(J_A=j) \leq  16 \,  \sum_{1 \leq j<n}  \frac 1 {j^2}  \, \frac{ j }{ (n-1)}\cdot \frac 1 { \binom n i} \ll \frac 1 { \binom n i}  \cdot \frac {\log n} n. \notag
\end{align*}
This leads to
\begin{align}\label{M_first_term_4}
\sum_{A}\VV\Big(\1_{\mathscr F_A }\Big) \cdot \EE( L_A^2) & \ll \sum_{r< i\leq a_{n,h}} \binom n i \cdot i \Big(\frac h {h+2}\Big)^{i-1} \cdot \frac 1 { \binom n i}  \cdot \frac {\log n} n \leq \varepsilon_1(r) \,  \frac {\log n} n 
\end{align}
where $ \varepsilon_1(r) := \sum_{i=r+1}^{\infty} i \Big(\frac h {h+2}\Big)^{i-1} $ has the property that
\be\label{epsilon_1}
\varepsilon_1(r)= o(1) \qquad \text{ when } r \to \infty.
\ee

Putting together (\ref{M_first_term_1}) and (\ref{M_first_term_2}) - (\ref{M_first_term_4}) we obtain for the first term on the right-hand side of (\ref{variance_M}) that
\begin{align}\label{M_first_term_conclusion}
\sum_{ r<|A|, |A'|\leq a_{n,h}}\EE&\Big((\1_{\mathscr F_A } - \EE(\1_{\mathscr F_A }))(\1_{\mathscr F_{A' }} - \EE(\1_{\mathscr F_{A'} }))\Big) \cdot \EE( L_A L_{A'})  \ll \Big(\varepsilon_1(r) \,  \frac {\log n} n + \frac 1 n \Big).
\end{align}

In order to obtain the claim of the Proposition we are left to bound the second term on the right-hand side of (\ref{variance_M}). For this term it holds that
\begin{align}\label{M_second_term_1}
\VV\Big(\sum_{r<i\leq a_{n,h}} p_{i,h}^n \cdot \ML^{n,i}\Big)  = \sum_{r<i\leq a_{n,h}} &(p_{i,h}^n)^2 \cdot \VV\Big(\ML^{n,i}\Big) + \sum_{r<i, i'\leq a_{n,h} \atop i \neq i'} p_{i,h}^n\,p_{i',h}^n \cdot \mathbb{COV}\Big(\ML^{n,i}, \ML^{n, i'}\Big).
\end{align}

The variances and covariances of the internal lengths of different orders can be easily obtained from the results of Fu \cite{Fu95} on the variances and covariances of the numbers $M_i(n)$ of mutations carried by exactly $i$ individuals in a population of size $n$ evolving according to the Moran model under the infinitely many sites mutation model. In this setting mutations are modelled as points of a Poisson process with constant rate $\frac \phi 2$ per unit length on the branches of the coalescent tree. Therefore it holds that
\[\Big(M_i(n) \mid \ML^{n, i}\Big) \sim \text{ Poisson }\Big(\frac \phi 2 \,\ML^{n, i}\Big).\]
We thus obtain by the law of total variance and of total covariance respectively that 
\begin{align}\label{var_L}
\VV(M_i(n)) & = \VV\Big(\EE(M_i(n) \mid \ML^{n, i})\Big)+ \EE\Big(\VV(M_i(n) \mid \ML^{n, i})\Big) \notag
\\ 
& = \VV\Big(\frac \phi 2 \,\ML^{n, i}\Big)+\EE\Big(\frac \phi 2 \,\ML^{n, i}\Big) \notag
\\ 
& =\frac {\phi^2} 4 \VV\Big(\ML^{n, i}\Big) + \frac \phi i
\end{align}
and due to the independence ensured by the Poisson structure
\begin{align}\label{covar_L}
\mathbb{COV}(M_i(n), M_{i'}(n)) & =\mathbb{COV}\Big(\EE(M_i(n) \mid \ML^{n, i},\ML^{n, i'}), \EE(M_{i'}(n) \mid \ML^{n, i}, \ML^{n, i'})\Big) \notag
\\ 
& \hphantom{1234567890}+ \EE\Big(\mathbb{COV}(M_i(n), M_{i'}(n) \mid \ML^{n, i},\ML^{n, i'})\Big) \notag
\\ 
& =\mathbb{COV}\Big(\frac \phi 2 \,\ML^{n, i}, \frac \phi 2 \,\ML^{n, i'}\Big)+ 0 \notag
 \\ 
&=\frac {\phi^2} 4 \mathbb{COV}\Big(\ML^{n, i}, \ML^{n, i'}\Big). 
\end{align}

The results of \cite{Fu95} say that
\be\label{var_covar_Fu}
\VV(M_i(n))= \phi^2 \sigma_{ii} + \frac \phi i\qquad\text{ and }\qquad\mathbb{COV}(M_i(n), M_{i'}(n))= \phi^2 \sigma_{ii'},
\ee
where in particular for $ i <\frac n 2$
\[
\sigma_{ii} =\beta_n(i+1)
\]
and for $i>i'$, $ i +i'<\frac n 2$
\[
\sigma_{ii'} = \frac{\beta_n(i+1)-\beta_n(i)} 2
\]
with $h_n=\sum_{j=1}^{n-1} \frac 1 j$ and
\[\beta_n(i) = \frac {2n}{(n-i+1)(n-i)}(h_{n+1}-h_i)- \frac 2 {n-i}.\]
Therefore from (\ref{var_L}), (\ref{covar_L}) and (\ref{var_covar_Fu}) it follows that 
\[\VV\Big(\ML^{n, i}\Big) = 4 \sigma_{ii} \qquad\text{ and }\qquad\mathbb{COV}\Big(\ML^{n, i}, \ML^{n, i'}\Big)= 4 \sigma_{ii'}.\]
Turning now to (\ref{M_second_term_1}) note from the definition of $a_{n,h}$ given in (\ref{def_a_nh}) that the indices $i$ and $i'$ run only up to values of the order $\log n$. Therefore, for large $n$, it holds that $i < \frac n 2$ and $i+i' < \frac n 2$ and hence we need the values of $\sigma_{ii}$ and $\sigma_{ii'}$ only in these particular cases. For $i < \frac n 2$ the following bound holds
\[\sigma_{ii} = \beta_n(i+1) \leq \frac {2n\,h_{n+1}}{(n-i+1)(n-i)} \leq c \cdot \frac {\log n} {n},\]
for $c$ a finite constant independent of $i$. Also $ \sigma_{ii'}\leq 0$ for $i>i'$, $i+i' < \frac n 2$ (see (36) in \cite{Fu95}). Hence using (\ref{expectation_1_SA}), (\ref{M_second_term_1}) becomes
\begin{align}\label{M_second_term_conclusion}
\VV\Big(\sum_{r<i\leq a_{n,h}} p_{i,h}^n \cdot \ML^{n,i}\Big) &\leq  \sum_{r<i\leq a_{n,h}} c \cdot 4i^2\Big(\frac h {h+2}\Big)^{2i-2} \cdot \frac {\log n} {n} \leq \varepsilon_2(r) \,  \frac {\log n} n,
\end{align}
where $ \varepsilon_2(r) := c \sum_{i=r+1}^{\infty} 4i^2 \Big(\frac h {h+2}\Big)^{2i-2} $ has the property that
\be\label{epsilon_2}
\varepsilon_2(r)= o(1) \qquad \text{ when } r \to \infty.
\ee

Putting now together (\ref{variance_M}), (\ref{epsilon_1}), (\ref{M_first_term_conclusion}), (\ref{M_second_term_conclusion}) and (\ref{epsilon_2}) we obtain the claim.
\ep

\section{Proof of Theorem 1}\label{Proof of the main result}

\subsection{Preliminaries}\label{subI}

\bprop
For any $r \in \NN$ it holds that
\be\label{L_A_to_to_2_J_A_with_F}
\sum_{1 \leq  |A| \leq r} (\1_{\mathscr F_A} \cdot L_A - \EE(\1_{\mathscr F_A} \cdot L_A)) = \sum_{1 \leq |A|\leq r}  \1_{\mathscr F_A}\Big(\frac 2 {J_A} - \EE\Big(\frac 2 {J_A}\Big)\Big)+ O_P(n^{-\nicefrac 1 2})
\ee
and also
\begin{equation}\label{L_A_to_to_2_J_A_without_F}
\sum_{1 \leq  |A| \leq r} (L_A - \EE(L_A)) = \sum_{1 \leq |A|\leq r} \Big(\frac 2 {J_A} - \EE\Big(\frac 2 {J_A}\Big)\Big)+ O_P(n^{-\nicefrac 1 2}).
\ee
\eprop

\bp
We start by making the observation that for any $r \in \NN$ the following claim holds:
\be\label{th_1_obs1}
\sum_{1 \leq  |A| \leq r} (\1_{\mathscr F_A} \cdot L_A - \EE(\1_{\mathscr F_A} \cdot L_A)) = \sum_{1 \leq  |A| \leq r}  \1_{\mathscr F_A} \cdot (L_A- \EE(L_A)) + O_P(n^{-\nicefrac 1 2}).
\ee
Indeed, since the collections of random variables $\{\1_{\mathscr F_A}\}_{A \in \MP_n}$ and $\{L_A\}_{A \in \MP_n}$ are independent we have that
\begin{align*}
\sum_{1 \leq  |A| \leq r} (\1_{\mathscr F_A}& \cdot L_A - \EE(\1_{\mathscr F_A} \cdot L_A)) =  \sum_{1 \leq  |A| \leq r}  \1_{\mathscr F_A} \cdot (L_A- \EE(L_A))  + \sum_{1 \leq  |A| \leq r}  (\1_{\mathscr F_A} - \EE(\1_{\mathscr F_A})) \cdot \EE(L_A).
\end{align*}
Using \eqref{bound_covariances} and \eqref{expectation_L_A} we obtain
\begin{align*}
\VV\Big(\sum_{1 \leq  |A| \leq r}  (\1_{\mathscr F_A} - \EE(\1_{\mathscr F_A})) \cdot \EE(L_A)\Big) & = \sum_{1 \leq |A|, |A'| \leq r} \mathbb{COV}(\1_{\mathscr F_A} ,\1_{\mathscr F_{A'}}) \cdot \EE(L_A)\EE(L_{A'}) 
\\ 
& \ll  \sum_{i, i'=1}^r  \sum_{ |A|=i, |A'|=i'\atop A\cap A' =\emptyset }  \frac 1 n \cdot \frac {1} {\binom n i \binom n {i'}}  +  \sum_{i, i'=1}^r  \sum_{|A|=i, |A'|=i'\atop A\cap A' \neq \emptyset } \frac {1} {\binom n i \binom n {i'}}.
\end{align*}
We now use the fact that the number of non-disjoint sets with cardinalities $i$ and $i'$ is 
\begin{align}\label{no_non_disjoint_sets}
\Big| \{A, A' \in \MP_n, |A|=i, & |A'|=i \text{ such that } A\cap A' \neq \emptyset\}\Big| \leq \binom n i \cdot i \binom n {i'-1} \ll \frac{1}{n} \binom n i \binom n {i'}
\end{align}
and obtain that 
\begin{align*}
\VV\Big(\sum_{1 \leq  |A| \leq r}  (\1_{\mathscr F_A} - \EE(\1_{\mathscr F_A})) \cdot \EE(L_A)\Big) & \ll \frac 1 n.
\end{align*}
Thus \eqref{th_1_obs1} is true.

\vspace{0.8cm}
Let us now introduce some more notation. For $i \in \NN$ let 
\[\mathcal G_i:=\sigma\Big(J_A, K_A : A \in \MP_n, |A|=i\Big)\]
be the $\sigma$-algebra containing the information about the final and initial levels of the branches of order $i$ in the coalescent at time 0 and let
\be\label{def_G}
\mathcal G:=\sigma\Big(\bigcup_{i=1}^{n-1} \mathcal G_i\Big).
\ee
The $\sigma$-algebra $\mathcal G$ contains the whole information on the structure of the coalescent tree (up to the exponential times). The conditional expectation of the length of the branch supporting the leaves with labels in $A$ is
\be\label{L_A_cond_G}
\EE(L_A \mid \mathcal G)=\EE\Big(\sum_{j=J_A+1}^{K_A} X_j\mid \mathcal G\Big)=\sum_{j=J_A+1}^{K_A} \EE(X_j) =  \frac 2 {J_A} - \frac 2 {K_A}.
\ee

We next show that by replacing the exponential times $X_j$ in the lengths $L_A$ appearing on the right-hand side of \eqref{th_1_obs1} by their expectations leads to a negligible error as the total population size $n$ tends to infinity. In other words, the randomness brought in by the inter-coalescent times can be neglected for big population sizes. The information on the tree structure is contained in $\mathcal G$ and therefore the observation we just made amounts to saying that for $i \geq 1$
\begin{equation}\label{L_A_to_cond_expectation}
\sum_{ |A|=i}  \1_{\mathscr F_A} L_A=\sum_{|A|=i}  \1_{\mathscr F_A} \EE(L_A \mid \mathcal G) + O_P(n^{-\nicefrac 1 2}).
\end{equation}

Indeed since the collection of random variables $\{ \1_{\mathscr F_A}\}_{A\in \MP_n}$ is independent of $\{L_A\}_{A\in \MP_n}$ and independent of the $\sigma$-algebra $\mathcal G$, it holds that 
\begin{align*}
\EE\Big(\Big(\sum_{ |A|=i}&  \1_{\mathscr F_A}(L_A-\EE(L_A \mid \mathcal G))\Big)^2 \mid \mathcal G\Big)
\\
& = \EE\Big( \sum_{ |A|=|A'|=i}  \1_{\mathscr F_A}  \1_{\mathscr F_{A'}}(L_A-\EE(L_A \mid \mathcal G))(L_{A'}-\EE(L_{A'} \mid \mathcal G)) \mid \mathcal G\Big)
\\
& =  \sum_{ |A|=|A'|=i} \EE(\1_{\mathscr F_A}\1_{\mathscr F_{A'}})\mathbb{COV}(L_A, L_{A'} \mid \mathcal G).
\end{align*}
Note that if $A \cap A' \notin \{A, A', \emptyset\}$ then there cannot be two branches in the coalescent, one supporting the leaves with labels in $A$ and one the leaves with labels in $A'$ and therefore in such case $ L_A$ or $L_{A'}$ is (by definition) equal to 0 and thus $\EE(L_A L_{A'} \mid \mathcal G)=0$ and $\mathbb{COV}(L_A, L_{A'} \mid \mathcal G)\leq0$. Using now \eqref{def_SA} we obtain that
\begin{align*}
\EE\Big(\Big(\sum_{ |A|=i} & \1_{\mathscr F_A}(L_A-\EE(L_A \mid \mathcal G))\Big)^2 \mid \mathcal G\Big) 
\\
& \leq  \sum_{ A } \VV\Big(\sum_{j=J_A+1}^{K_A}X_j \mid \mathcal G\Big)+ \sum_{ A \cap A' = \emptyset} \Big ( \mathbb{COV}\Big(\sum_{j=J_A+1}^{K_A}X_j , \sum_{j'=J_{A'}+1}^{K_{A'}}X_{j'}\mid \mathcal G\Big) \Big )^+ 
\\
& \leq \sum_{A } \sum_{j=J_A+1}^{K_A}\VV(X_j) + \sum_{ A \cap A' = \emptyset} \sum_{j=J_A+1}^{K_A} \sum_{j'=J_{A'}+1}^{K_{A'}}\Big (  \mathbb{COV}(X_j, X_{j'})\Big )^+.
\end{align*}
Recall that the inter-coalescent times $X_j$ are independent and exponentially distributed with parameter $\binom j 2$. Therefore
\begin{align*}
\EE\Big(\Big(\sum_{A} \1_{\mathscr F_A}(L_A-\EE(L_A \mid \mathcal G))\Big)^2 \mid \mathcal G\Big) &\leq \sum_{ A } \sum_{j=J_A+1}^{K_A}\VV(X_j) +  \sum_{ A \cap A' = \emptyset} \sum_{j=J_A \vee J_{A'}+1}^{K_A \wedge K_{A'}} \VV(X_j)
\\ 
& = \sum_{ A} \sum_{j=J_A+1}^{K_A}\frac 1 {\binom j 2 ^2}+  \sum_{ A \cap A' = \emptyset} \sum_{j=J_A \vee J_{A'}+1}^{K_A \wedge K_{A'}} \frac 1 {\binom j 2 ^2}
\\ 
& \ll \sum_{ A } \frac 1{J_A^3}+ \sum_{  A \cap A' = \emptyset}\frac 1 {J_A^3 \vee J_{A'}^3}.
\end{align*}

Taking the expectation and using Lemma \ref{lemma_J} we obtain that
\begin{align*}
\EE\Big(\Big(\sum_{ |A|=i}  \1_{\mathscr F_A}(L_A-\EE(L_A \mid \mathcal G))\Big)^2\Big)& \leq \sum_{ A} \sum_{1 \leq j <n} \frac 1 {j^3} \cdot \frac{2 j }{(n-1)} \cdot \frac 1 {\binom n {i}}
\\
&  \quad \quad + \sum_{  A \cap A'= \emptyset} \sum_{1 \leq j <n} \sum_{j'=1}^{j} \frac 1 {j^3} \cdot \frac{4 jj' }{(n-1)(n-2)} \cdot \frac 1 {\binom n {i, i, n-2i}} 
\\
& \quad \quad +\sum_{  A \cap A'= \emptyset} \sum_{1 \leq j <n}\sum_{j'\geq j+1} \frac 1 {j'^3} \cdot \frac{4 jj' }{(n-1)(n-2)} \cdot \frac 1 {\binom n {i, i, n-2i}}
\\
& \ll \frac{1 }{n}.
\end{align*}
This yields \eqref{L_A_to_cond_expectation}.

Recall from \eqref{L_A_cond_G} that $\EE(L_A \mid \mathcal G)=\frac 2 {J_A} - \frac 2 {K_A}$.
In the next step we show that the terms $\frac 2 {K_A}$, $A \in \MP_n$, from this representation can be neglected when investigating the evolution of the external length. Using \eqref{L_A_to_cond_expectation}, this amounts to saying that for $i \geq 1$ the following claim holds:
\begin{equation}\label{L_A_to_2_J_A}
\sum_{|A|=i}  \1_{\mathscr F_A} (L_A- \EE(L_A))=\sum_{|A|=i}  \1_{\mathscr F_A}\Big(\frac 2 {J_A} - \EE\Big(\frac 2 {J_A}\Big)\Big)+ O_P(n^{-\nicefrac 1 2}).
\end{equation}
In order to see this it suffices to show that
\be\label{without_K_A}
\EE\Big(\Big(\sum_{|A|=i}  \1_{\mathscr F_A} \Big(\frac 2 {K_A} - \EE\Big(\frac 2 {K_A}\Big)\Big)\Big)^2\Big) \ll \frac 1 n.
\ee
The claim holds for $i=1$ since by definition $K_A=n$ for all sets with cardinality one. In the case $i\geq 2$ remember that  if $A \cap A' \notin \{A, A', \emptyset\}$ then there cannot be two branches in the coalescent, one supporting the leaves with labels in $A$ and one the leaves with labels in $A'$ and therefore in such case $K_A$ or $K_{A'}$ is (by definition) equal to $\infty$ and thus $\mathbb{COV}(\frac 2 {K_A}, \frac 2 {K_{A'}})\leq 0$. 

Using Lemma \ref{lemma_J} we obtain for $i\geq 2$ that
\begin{align*}
\EE\Big(\Big(\sum_{|A|=i}  &\1_{\mathscr F_A} \Big(\frac 2 {K_A} - \EE\Big(\frac 2 {K_A}\Big)\Big)\Big)^2\Big) 
\\ 
&\leq \sum_{|A|=i} \VV\Big(\frac 2 {K_A}\Big)+  \sum_{|A|, |A'|=i \atop A \cap A' =\emptyset} \EE( \1_{\mathscr F_A}  \1_{\mathscr F_{A'}})\cdot \mathbb{COV}\Big(\frac 2 {K_A}, \frac 2 {K_{A'}}\Big)
\\
& \ll \sum_{|A|=i}\sum_{k=2}^n \frac 1 {k^2}\PP(K_A=k)
\\
&\quad\quad+\sum_{A \cap A' =\emptyset} \EE( \1_{\mathscr F_A}  \1_{\mathscr F_{A'}}) \sum_{k= 2}^n \sum_{k'= 2}^n \frac 1 {kk'} \big( \PP(K_A=k, K_{A'}=k')- \PP(K_A=k) \PP(K_{A'}=k')\big)
\\
& \ll \sum_{|A|=i}\sum_{k=2}^n \frac 1 {k^2}\PP(K_A=k) +\sum_{A \cap A' =\emptyset}\sum_{k= 2}^n  \sum_{k'= 2}^n\frac 1 {kk'} \PP(K_A=k) \PP(K_{A'}=k')\cdot \frac 1 n
\\
& \ll \sum_{|A|=i}\sum_{k=2}^n \frac 1 {k^2}\cdot \frac{k^2}{(n-1)(n-2)}\cdot \frac 1 {\binom n i} + \sum_{A \cap A' =\emptyset}\sum_{k= 2}^n  \sum_{k'= 2}^n \frac 1 {kk'}\cdot\frac {k^2k'^2}{(n-1)^2(n-2)^2}\cdot\frac 1 {{\binom n i}^2}\cdot\frac 1 n 
\\
& \ll  \frac 1 n
\end{align*}
and thus \eqref{without_K_A} holds. Putting now together \eqref{th_1_obs1} and \eqref{L_A_to_2_J_A} we obtain the first claim of the proposition.

The calculations above remain valid if we replace $\1_{\mathscr F_A}$ by 1. Thus, in the same manner, we also obtain the second claim. 
\ep

\subsection{From dependent to independent coefficients}\label{dep_to_indep_thinning}

This subsection contains a key building block of the proof of Theorem \ref{mainresult_fdd}, namely we show that the dependent random coefficients $\{\1_{\mathscr F_A}\}_{A \in \MP_n}$ which appear in the sum on the right-hand side of \eqref{L_A_to_to_2_J_A_with_F} can be replaced in the infinite population size limit by random coefficients that are independent of one another and of the coalescent at time 0. 

We fix $r \in \NN$ and group the branches of the coalescent at time 0 which support at most $r$ leaves (and thus the sets $A \in \MP_n, 1 \leq |A|\leq r$) according to whether their initial and final levels are above or below level $ b_n$ with
\[b_n:= \Big\lceil \frac n {\log n}\Big\rceil.\]
Recall that if there exists no branch in the coalescent at time 0 supporting the leaves with labels in a set $A$ then the initial and final levels $K_A$ and $J_A$ are both equal to $\infty$.

Recall also that $\1_{\mathscr F_A}=\1_{\{Z^A_h=1\}}$. If the branch supporting the leaves with labels in the set $A$ with $|A|=i$ is formed above level $b_n$, i.e. $K_A < b_n$, or ends below or at level $b_n$, i.e. $J_A \geq b_n$, we replace the coefficient $\1_{\mathscr F_A}$ by $\widetilde W^A_{h}:=\1_{\{\widetilde B^A_h=1\}}$, where the processes $\widetilde B^A$ are independent birth and death processes started with $i$ individuals at time 0 with birth and death rates equal to $\frac 1 2 k$ when the process is in state $k$. Moreover, we choose the birth and death processes to be independent of one another and independent of the coalescent at time 0 and its evolution forwards in time. 

Let us now look at the branches that are formed below level $b_n$ and end above this level. These are the $b_n$ branches extant at level $b_n$. Observe that they support leaves with labels in {\it disjoint} sets. Out of these random sets, let $A_1, \dots, A_{l_n}$, $l_n \leq b_n$, denote those which have cardinalities less than or equal to $r$: $|A_l|=:i_l\leq r$, $1\leq l\leq l_n$. Let 
\[
\bm\bar A:=\bigcup_{l=1}^{l_n} A_j \qquad \text{ and } \qquad \bm\bar i:= i_1+\dots+i_{l_n}
\]
and note that
\be\label{A_and_i_bar}
\bm\bar i \leq r \cdot \Big\lceil \frac n {\log n}\Big\rceil.
\ee

Now we describe our coupling construction. Given the $\sigma$-algebra $\mathcal G$, we couple, like in Section \ref{A birth and death process}, independent birth and death processes $\{B^{j}\}_{j \in \bm\bar A}$ to the processes $\{Z^{j}\}_{j \in \bm\bar A}$. In particular, the processes $\{B^{A_l}\}_{1 \leq l\leq l_n}$ are coupled to the processes  $\{Z^{A_l}\}_{1 \leq l\leq l_n}$ (the random time change depending on the set $\bm\bar A$).  This gives rise to the random variables $W^{A_l}_{h}:= \1_{\{B^{A_l}_h=1\}}$ which we use to replace the coefficients $\1_{\mathscr F_{A_l}}$. Since the sets $A_1, \dots, A_{l_n}$ are all disjoint, it follows by Lemma \ref{indep_B_j} that the processes $B^{A_1}, \dots, B^{A_{l_n}}$, and thus the random variables $W^{A_1}_{h}, \dots, W^{A_{l_n}}_{h}$, are independent of one another and independent of the coalescent at time 0.

Through this procedure we obtain the collection of random coefficients
\be\label{def_V_A}
V^A_{h}:=  \1_{\{J_A < b_n\leq K_A\}}\cdot W^A_{h}+\1_{\{K_A < b_n\}\cup \{J_A \geq b_n\}} \cdot  \widetilde W^A_{h}, \qquad\qquad A \in \mathcal P_{n,r},
\ee
where $\mathcal P_{n,r}:=\{A \subset \{1, \dots, n\}, 1 \leq |A|\leq r\}$. Note that from construction the random variables $V^A_{h}$ are Bernoulli$(p_{i,h})$-distributed if the set $A$ has cardinality $i$ with $p_{i,h}:=\PP(B_h^{A}=1 \mid B_0^{A}=i)$. These are the coefficients we use in order to replace the $\1_{\mathscr F_A}$'s. In what follows we show that they are independent from one another and independent of the coalescent tree at time 0 (even though $J_A$ and $K_A$ are used in their definition). 

Since the dependence on the coalescent tree at time 0 can come only from the tree structure (and not from the exponential times) and the whole information on the tree structure is contained in the $\sigma$-algebra $\mathcal G$ defined in \eqref{def_G}, for the latter property we only need to prove the independence of the collection $\{V^A_{h}\}_{A \in \mathcal P_{n,r}}$ of the $\sigma$-algebra $\mathcal G$. To this aim consider bounded, measurable functions $\{f_A\}_{A \in \mathcal P_{n,r}}$ depending on the sets $A$. Then
\begin{align*}
\EE\Big[&\prod_{A \in \mathcal P_{n,r}} f_A(V^A_h) \mid \mathcal G\Big]
\\
&= \EE\Big[\prod_{A \in \mathcal P_{n,r}} \Big( \1_{\{J_A < b_n\leq K_A\}}\cdot f_A(W^A_h)+\1_{\{K_A < b_n\}\cup \{J_A \geq b_n\}} \cdot  f_A(\widetilde W^A_h)\Big) \mid \mathcal G\Big] \notag
\\
&= \EE\Big[\sum_{\mathcal A \in \mathcal P (\mathcal P_{n,r})} \Big(\prod_{A \in \mathcal A} \big( \1_{\{J_A < b_n\leq K_A\}}\cdot f_A(W^A_h)\big) \cdot \prod_{A \notin \mathcal A} \big( \1_{\{K_A < b_n\}\cup \{J_A \geq b_n\}} \cdot  f_A(\widetilde W^A_h)\big)\Big) \mid \mathcal G\Big]\notag
\\
&=\sum_{\mathcal A \in \mathcal P (\mathcal P_{n,r})}\Big( \prod_{A \in \mathcal A}  \1_{\{J_A < b_n\leq K_A\}}  \cdot \prod_{A \notin \mathcal A} \1_{\{K_A < b_n\}\cup \{J_A \geq b_n\}} \cdot\EE\Big[ \prod_{A \in \mathcal A} f_A(W^A_h) \cdot \prod_{A \notin \mathcal A} f_A(\widetilde W^A_h) \mid \mathcal G\Big] \Big),
\end{align*} 
where the last equality follows from the fact that the indicator functions are $\mathcal G$-measurable. 

Recall now that if there are two non-disjoint sets $A_1, A_2 \in \mathcal A$, then $\1_{\{J_{A_1} < b_n\leq K_{A_1}\}}  \cdot \1_{\{J_{A_2} < b_n\leq K_{A_2}\}} =0$ and therefore the first product in the parenthesis on the right-hand side of the above equality is equal to 0. Using this observation together with the independence of the random variables $W^A_h$ for disjoint sets $A$, the independence of the random variables $\widetilde W^A_h$ and their independence of the $\sigma$-algebra $\mathcal G$ we obtain that
\begin{align*}
\EE\Big[\prod_{A \in \mathcal P_{n,r}} f_A(V^A_h) \mid \mathcal G\Big] &=\sum_{\mathcal A \in \mathcal P (\mathcal P_{n,r}) \atop \mathcal A \text{ contains only disjoint sets}}  \prod_{A \in \mathcal A}  \1_{\{J_A < b_n\leq K_A\}} \cdot \prod_{A \in \mathcal A} \EE[ f_A(W^A_h) ]
\\
&\qquad\qquad\qquad\qquad\qquad\cdot\prod_{A \notin \mathcal A} \1_{\{K_A < b_n\}\cup \{J_A \geq b_n\}} \cdot\prod_{A \notin \mathcal A} \EE[ f_A(\widetilde W^A_h)].
\end{align*}  
Now we can reverse these steps and get that
\begin{align*}
\EE\Big[\prod_{A \in \mathcal P_{n,r}} f_A(V^A_h) \mid \mathcal G\Big]&= \prod_{A \in \mathcal P_{n,r}} \EE\Big[ \Big( \1_{\{J_A < b_n\leq K_A\}}\cdot f_A(W^A_h)+\1_{\{K_A < b_n\}\cup \{J_A \geq b_n\}} \cdot  f_A(\widetilde W^A_h)\Big)\Big]
\\
&= \prod_{A \in \mathcal P_{n,r}} \EE[ f_A(V^A_h)]
\end{align*}
which proves the independence of the collection $\{V^A_h\}_{A \in \mathcal P_{n,r}}$ of the $\sigma$-algebra $\mathcal G$ and thus of the coalescent tree at time 0. Taking now the expectation in the equality above we also obtain the independence of the random variables $V^A_h$ from one another.

\vspace{0.8cm}
We next show that in the infinite population size limit one can replace the coefficients $\1_{\mathscr F_A}$  by the random variables $V^A_h$, namely that the following holds:
\begin{align}\label{dep_to_indep}
 \sum_{ |A| \leq r} & \1_{\mathscr F_A} \cdot  \Big(\frac 2 {J_A} - \EE\Big(\frac 2 {J_A}\Big)\Big) =  \sum_{ |A|\leq r}V^A_h \cdot  \Big(\frac 2 {J_A} - \EE\Big(\frac 2 {J_A}\Big)\Big)+ O_P\Big( \Big( \frac {\log \log n}{n}\Big)^{\nicefrac 1 2}\Big).
\end{align}

We consider the three groups of branches (for which $K_A < b_n$, $J_A \geq b_n$ and respectively $J_A < b_n\leq K_A$) separately. For each $1 \leq i\leq r$ fixed we have due to independence that
\begin{align*}
 \EE\Big( \Big(\sum_{ |A|=i} &(\1_{\mathscr F_A} - V^A_h) \cdot \Big(\frac 2 {J_A} - \EE\Big(\frac 2 {J_A}\Big)\Big) \cdot  \1_{\{K_A < b_n\}}\Big)^2 \Big)  \notag
\\
& = \sum_{ |A|=|A'|=i} \EE \Big((\1_{\mathscr F_A} -V^A_h) (\1_{\mathscr F_{A'}} - V_h^{A'})\Big)  \notag
\\
& \quad \quad \cdot \EE \Big( \Big(\frac 2 {J_A} - \EE\Big(\frac 2 {J_A}\Big)\Big)  \Big(\frac 2 {J_{A'}} - \EE\Big(\frac 2 {J_{A'}}\Big)\Big) \Big) ; K_A, K_{A'} < b_n \Big). \notag
\end{align*}
Remember that if $A\cap A' \notin \{A, A', \emptyset\}$ then at least one of the branches with leaves in $A$, respectively $A'$, does not exist and by definition $K_A=\infty$ or $K_{A'}=\infty$. Therefore the event $\{K_A, K_{A'} <b_n\}$ is empty and thus in the sum above the summands for which $A\cap A' \notin \{A, A', \emptyset\}$ are equal to 0.  We obtain that
\begin{align}\label{group1_intermed}
 \EE\Big( &\Big(\sum_{ |A|=i} (\1_{\mathscr F_A} - V^A_h) \cdot \Big(\frac 2 {J_A} - \EE\Big(\frac 2 {J_A}\Big)\Big)  \cdot  \1_{\{K_A < b_n\}}\Big)^2 \Big)  
  \\
 & \leq \sum_{ A} \EE \Big(\frac 1 {J_A^2}+\EE\Big(\frac 1 {J_A}\Big)^2\Big); K_A < b_n \Big) + \sum_{ A\cap A' = \emptyset } \Big(\mathbb{COV} \Big(\frac 1 {J_A}\cdot \1_{\{K_A < b_n\}}, \frac 1 {J_{A'}}\cdot \1_{\{K_{A'} < b_n\}}\Big)\Big)^+.\notag
\end{align}

Without loss of generality we can choose $A:=\{1, \dots, i\}$ and $A':=\{i+1, \dots, 2i\}$. Let us write for short $K:=K_A$, $J:=J_A$,  $K':=K_{A'}$ and $J':=J_{A'}$. 

Note that $K_A<b_n$, respectively  $K_{A'}<b_n$ implies that $i \geq 2$, respectively $i'\geq 2$, and thus from Lemma \ref{lemma_J}
\begin{align}\label{cov_J,J'}
&\mathbb{COV}\Big(\frac 1 {J_A}\cdot \1_{\{K_A < b_n\}}, \frac 1 {J_{A'}}\cdot \1_{\{K_{A'} < b_n\}}\Big) \notag
 \\
 &= \sum_{k, k' \leq b_n} \sum_{1\leq j\leq k-1 \atop 1\leq j'\leq k'-1} \frac 1 {jj'} \cdot \big(\PP(K=k, J=j, K'=k', J'=j')- \PP(K=k, J=j)\PP(K'=k', J'=j')\big) \notag
 \\ 
 & \leq  c \cdot \frac 1 n \sum_{k, k' \leq b_n} \sum_{1\leq j\leq k-1 \atop 1\leq j'\leq k'-1} \frac 1 {jj'} \cdot\PP(K=k, J=j)\PP(K'=k', J'=j') \notag
  \\ 
 & \leq  c \cdot \frac 1 n \sum_{k, k' \leq b_n} \sum_{1\leq j\leq k-1 \atop 1\leq j'\leq k'-1} \frac 1 {jj'} \cdot \frac {jj'}{(n-1)^2(n-2)^2}\cdot \frac 1 {{\binom n i}^2} \notag
 \\ 
 & \leq  c \cdot \frac 1 n \cdot \frac 1 {{\binom n i}^2} 
\end{align}
for $c<\infty$ depending on $i$ and $i'$.

Therefore using \eqref{expect_1/J} and again Lemma \ref{lemma_J} the inequality \eqref{group1_intermed} becomes
\begin{align}\label{var_K_less_bn}
\EE&\Big( \Big(\sum_{ |A|=i} (\1_{\mathscr F_A} - V^A_h) \cdot \Big(\frac 2 {J_A} - \EE\Big(\frac 2 {J_A}\Big)\Big) \cdot  \1_{\{K_A < b_n\}}\Big)^2 \Big)
\\
& \ll \binom n i  \sum_{k \leq b_n} \sum_{j=1}^{k-1} \Big(\frac 1 {j^2}+ \frac 1 {\binom n i ^2}\Big) \cdot \PP(K=k, J=j) + \binom n {i, i, n-2i} \cdot \frac 1 n \cdot \frac 1 {{\binom n i}^2} \notag
\\
&\ll \binom n i \sum_{k\leq b_n} \sum_{j=1}^{k-1}\frac 1 {j^2} \cdot \frac {2ij} {(n-1)(n-2)} \cdot \frac 1 {\binom n i} + \frac 1 n \notag
\\
&\ll \sum_{k\leq b_n} \frac {i\log k}{ n^2} +  \frac 1 n \notag
\\
&\ll \frac {ib_n \log b_n}{ n^2} +  \frac 1 n \notag
\\
&\ll \frac 1 n. \notag
\end{align}

For the branches for which the final level $J$ is greater than or equal to $b_n$ we use a similar argument. Again if $A\cap A' \notin \{A, A', \emptyset\}$ then by definition $J_A=\infty$ or $J_{A'}=\infty$ and thus  $\mathbb{COV}\Big(\frac 1 {J_A}\cdot \1_{\{J_A \geq b_n\}}, \frac 1 {J_{A'}}\cdot \1_{\{J_{A'} \geq b_n\}}\Big)\leq 0$. We obtain for $1 \leq i\leq r$ fixed that
\begin{align*}
\EE \Big( \Big(&\sum_{ |A|=i} (\1_{\mathscr F_A} - V^A_h) \cdot \Big(\frac 2 {J_A} - \EE\Big(\frac 2 {J_A}\Big)\Big)  \cdot  \1_{\{J_A \geq b_n\}} \Big)^2 \Big)
\\
& \ll \sum_{ A} \EE \Big(\frac 1 {J_A^2}+\EE\Big(\frac 1 {J_A}\Big)^2\Big); J_A \geq b_n \Big) + \sum_{ A\cap A' = \emptyset } \Big(\mathbb{COV}  \Big(\frac 1 {J_A}\cdot \1_{\{J_A \geq b_n\}}, \frac 1 {J_{A'}}\cdot \1_{\{J_{A'} \geq b_n\}}\Big)\Big)^+.
\end{align*}
As above, by Lemma \ref{lemma_J}
\begin{align}\label{cov_J,J'_2}
\mathbb{COV} \Big(\frac 1 {J_A}\cdot \1_{\{J_A \geq b_n\}}, \frac 1 {J_{A'}}\cdot \1_{\{J_{A'} \geq b_n\}}\Big)&= \sum_{b_n \leq j,j' <n} \frac 1 {jj'} \cdot \big(\PP(J=j, J'=j')- \PP(J=j)\PP(J'=j')\big) 
 \\ 
 & \ll \frac 1 n  \sum_{b_n \leq j,j' <n} \frac 1 {jj'} \cdot\PP(J=j)\PP(J'=j') \notag
  \\ 
 & \ll \frac 1 n  \sum_{b_n \leq j,j' <n} \frac 1 {jj'} \cdot \frac {jj'}{(n-1)^2}\cdot \frac 1 {{\binom n i}^2} \notag
 \\ 
 & \ll \frac 1 n \cdot \frac 1 {{\binom n i}^2} \notag
\end{align}
for $c<\infty$ depending on $i$ and $i'$ and thus
\begin{align}\label{var_J_greater_bn}
\EE\Big(\Big(  \sum_{|A|=i}  &(\1_{\mathscr F_A} - V^A_h) \cdot \Big(\frac 2 {J_A} - \EE\Big(\frac 2 {J_A}\Big)\Big) \cdot  \1_{\{J_A \geq b_n\}} \Big)^2 \Big) 
\\
&\ll \binom n i \sum_{b_n \leq j<n} \frac 1 {j^2} \cdot \frac {2j} {n-1} \cdot \frac 1 {\binom n i} + \frac 1 n  \notag
\\
&\ll \frac 1 n \cdot \log { \frac n { b_n}} +  \frac 1 n \notag
\\
&\ll \frac 1 n \cdot \log \log n. \notag
\end{align}

We are now left to consider the branches that are formed below and end above level $b_n$ in the coalescent. Again, using the same argument as above and the fact that for $A\cap A' \notin \{A, A', \emptyset\}$ the event $\{J_A < b_n\leq K_A, J_{A'} < b_n\leq K_{A'}\}$ is empty (since there can not be two branches in the tree that coexist at level $b_n$ and have common leaves), we have that
\begin{align}\label{var_between_1}
\EE\Big(& \Big(\sum_{ |A|=i}  (\1_{\mathscr F_A} - V^A_h) \cdot \Big(\frac 2 {J_A} - \EE\Big(\frac 2 {J_A}\Big)\Big) \cdot  \1_{\{ J_A < b_n\leq K_A\}} \Big)^2\Big) \notag
\\
& \ll \sum_{ A} \EE\Big((\1_{\mathscr F_A} - W^A_h)^2\Big) \cdot \EE \Big(\frac 1 {J_A^2}+ \EE\Big(\frac 1 {J_A}\Big)^2; J_A < b_n\leq K_A \Big) \notag
\\
& \quad \quad \quad \quad\quad \quad\qquad \quad + \sum_{ A\cap A' = \emptyset}\mathbb{COV} \Big(\frac 1 {J_A}\cdot \1_{\{J_A < b_n\leq K_A\}}, \frac 1 {J_{A'}}\cdot \1_{\{J_{A'} < b_n\leq K_{A'}\}}\Big)\notag
\\
& = \sum_{ A} \PP\Big(\1_{\{Z^A_h=1\}} \neq \1_{\{B^A_h=1\}}\Big) \cdot \EE \Big(\frac 1 {J_A^2}+ \EE\Big(\frac 1 {J_A}\Big)^2; J_A < b_n\leq K_A \Big)  \notag
\\
& \quad \quad \quad \quad\quad \quad\qquad \quad + \sum_{ A\cap A' = \emptyset}\mathbb{COV} \Big(\frac 1 {J_A}\cdot \1_{\{J_A < b_n\leq K_A\}}, \frac 1 {J_{A'}}\cdot \1_{\{J_{A'} < b_n\leq K_{A'}\}}\Big).
\end{align}
If the cardinality $i$ of the sets $A$ and $A'$ is at least 2, then the second sum on the right hand side can be bounded from above by $c \cdot \frac 1 n$ for a finite constant $c$ by using the same argument as the one used to obtain \eqref{cov_J,J'}. If $i=1$ recall that $K_A=n$  and thus $\{J_A < b_n\leq K_A\}=\{J_A < b_n\}$ (and the same for the set $A'$). By the argument used to obtain \eqref{cov_J,J'_2} we can again bound the sum from above by $c \cdot \frac 1 n$. As to the first sum on the right hand side of  \eqref{var_between_1} recall that we coupled all processes $Z^A$ with $A$ such that $J_A < b_n\leq K_A$ simultaneously and that the random time change depends just on the union $\bm\bar A$ of these sets.  From \eqref{prob_different} and \eqref{A_and_i_bar} we obtain that 
\[\PP\Big(\1_{\{Z^A_h=1\}} \neq \1_{\{B^A_h=1\}}\Big) \ll \frac {\bm\bar i} {n} \ll \frac 1 {\log n},\]
which plugged in \eqref{var_between_1} leads to 
\begin{align}\label{var_between_2}
\EE\Big(\Big(& \sum_{ |A|=i}  (\1_{\mathscr F_A} - V^A_h) \cdot\Big(\frac 2 {J_A} - \EE\Big(\frac 2 {J_A}\Big)\Big)\cdot  \1_{\{J_A < b_n\leq K_A\}} \Big)^2 \Big)
\\
& \ll \frac 1 {\log n} \sum_{A} \EE \Big(\frac 1 {J_A^2}+ \EE\Big(\frac 1 {J_A}\Big)^2 \Big) +  \frac 1 n \notag
\\
& \ll \frac 1 {\log n} \binom n i \sum_{j< n} \Big(\frac 1 {j^2}+ \frac 1 {\binom n i ^2}\Big) \cdot \PP(J=j) + \frac 1 n \notag
\\
&\ll \frac { \binom n i} {\log n} \sum_{j<  n} \frac 1 {j^2} \cdot \frac {2j} {(n-1)} \cdot \frac 1{\binom n i} +  \frac 1 n\notag
\\
&\ll \frac 1 {\log n} \cdot \frac {\log n}{ n}+ \frac 1 n \notag
\\
&\ll \frac 1 n. \notag
\end{align}
Putting now together \eqref{var_K_less_bn}, \eqref{var_J_greater_bn} and \eqref{var_between_2} we obtain \eqref{dep_to_indep}.

\subsection{Completion of the proof of the Theorem}\label{subIII}

Let $t>0$ be a fixed time measured on the evolutionary time-scale (the scale on which the lengths in the coalescent are measured) and recall that $h=\frac t n$. Remember from \eqref{ext_len_decomp} that, considering the coalescent tree of the population alive at time $h$, $I^n_{0,h}$ denotes the amount of external length in the coalescent at time $h$ that is gathered in the region near the leaves of the tree, more precisely in the time interval $(0,h)$. We start by showing that 
\be\label{I_to_show} 
I_{0,h}^n -\EE(I_{0,h}^n)=O_P(n^{-\nicefrac 1 2}).
\ee

Let $\widetilde X_j$, $2\leq j \leq n$, denote the inter-coalescence times in the coalescent of the population alive at time $h$ and consider the random time
\[T:= \sum_{j\geq \frac{2}{t+2}n}^n \widetilde X_j.\]
Since $\widetilde X_j \sim Exp\Big(\binom j 2\Big)$, it holds that $\EE(T)=\frac{t}{n}+ O\Big(\frac 1 {n^2}\Big)$ and $\VV(T) = O\Big( \frac{1}{n^3}\Big)$. Observe that the quantity $I^n_{0,\frac{t}{n}}-I^n_{0,T}$ is the external length gathered in an interval of length $|\frac{t}{n} -T|$. In this time interval some of the $n$ external branches may become internal. It follows that
\[| I^n_{0,\frac{t}{n}}-I^n_{0,T} | \leq n \cdot \Big |\frac{t}{n} -T\Big|\]
and thus 
\[ I^n_{0,\frac{t}{n}} - \EE(I^n_{0,\frac{t}{n}} ) = I^n_{0,T} -\EE(I^n_{0,T})+ O_P(n^{-\nicefrac 1 2}).\] 

In order to estimate $I^n_{0,T}$ we use Proposition 3 in \cite{JK11}. The arguments used in its proof apply completely unchanged also for $\alpha$ and $\beta$ depending on $n$. Choosing $\alpha$ such that $n^\alpha = \frac 2 {2+t} n$ and $\beta=1$ the proof gives
\[ \VV(I^n_{0,T}) = 8\Big(1- \frac{\log\big(\frac{2}{t+2} \, n\big)}{\log n}\Big)\cdot \frac{\log n}{n} + O(n^{-1})=O(n^{-1}),\]
and therefore \eqref{I_to_show} holds.

We are now ready for the final step. For ease of notation let us now consider just two times $0<h_1<h_2$. We show that for any $\alpha_1, \alpha_2 \in \RR$ the linear combination 
\[\sqrt{\frac {n}{4 \log n} }\Big(\alpha_1 (\ML_{h_1}^{n} - \EE(\ML_{h_1}^{n}))+ \alpha_2 (\ML_{h_2}^{n}- \EE(\ML_{h_2}^{n})\Big)\]
of the (centred and rescaled) lengths at the two times converges in distribution to a normal distributed random variable. The argument can be immediately extended for linear combinations of the lengths at times $0 < h_1< \dots< h_k$, $k \in \NN$.

Using the notation given at the beginning of Section \ref{Proof of the Proposition} we have that for each $r \in \NN$ 
\[\ML^{n}_{h_j}= \sum_{ |A| \leq r} \1_{\mathscr F_A(h_j)} \cdot L_A + \widetilde R^{n,r}_{h_j}+\bm\bar R^{n,r}_{h_j}+ I^n_{0,h_j}, \qquad j=1,2.\]
Here we write $\mathscr F_A(h_j)$ for the event that the number of descendants  at time $h_j$, $j \in \{1,2\}$, of the individuals with labels in the set $A$ at time 0 is equal to 1. The second and third terms have been treated in Proposition \ref{proposition_bigger_r} and we now focus on the first term. From \eqref{L_A_to_to_2_J_A_with_F} 
\begin{align*}
\sum_{ |A| \leq r} \1_{\mathscr F_A(h_j)} &\cdot L_A - \EE\Big(\sum_{ |A| \leq r} \1_{\mathscr F_A(h_j)} \cdot L_A \Big) = \sum_{ |A|\leq r}  \1_{\mathscr F_A(h_j)} \Big(\frac 2 {J_A} - \EE\Big(\frac 2 {J_A}\Big)\Big)+ O_P(n^{-\nicefrac 1 2}).
\end{align*}

Equation \eqref{dep_to_indep} allows to replace the dependent random coefficients $\1_{\mathscr F_A (h_1)}$ by the independent coefficients $V^A_{h_1}$ defined in \eqref{def_V_A}. In order to obtain the corresponding independent random coefficients for the time $h_2$ we let the birth and death processes $\widetilde W^A_{h_1}$ and $W^A_{h_1}$ used to define $V^A_{h_1}$ run further until time $h_2$ and define the coefficients $V^A_{h_2}$ as in \eqref{def_V_A} by using the states of the processes  $\widetilde W^A_{h_2}$ and $W^A_{h_2}$ at  time $h_2$. Note that the argument from the previous subsection used to prove \eqref{dep_to_indep} continues to hold also for the time $h_2$ and the coefficients defined as above. We thus obtain using \eqref{dep_to_indep}, \eqref{I_to_show} and Proposition \ref{proposition_bigger_r} that for $r \in \NN$
\begin{align*} 
\sqrt{\frac {n}{4 \log n} } &\Big(\alpha_1 (\ML_{h_1}^{n} - \EE(\ML_{h_1}^{n}))+ \alpha_2 (\ML_{h_2}^{n}- \EE(\ML_{h_2}^{n})\Big)
\\
&=\sqrt{\frac {n}{4 \log n} } \sum_{ |A| \leq r} \Big(\alpha_1 V^A_{h_1}+ \alpha_2 V^A_{h_2}\Big) \cdot \Big(\frac 2 {J_A} - \EE\Big(\frac 2 {J_A}\Big)\Big) \notag
\\
&\qquad+ \sqrt{\frac {n}{4 \log n} }\Big(\alpha_1\big(\widetilde R^{n,r}_{h_1} -\EE(\widetilde R^{n,r}_{h_1})\big)+ \alpha_2\big(\widetilde R^{n,r}_{h_2} -\EE(\widetilde R^{n,r}_{h_2})\big)\Big)\notag
\\
&\hphantom{123456789}+O_P\Big(\Big(\frac{\log\log n}{\log n}\Big)^{\nicefrac 1 2}\Big). \notag
\end{align*}

We now argue that the collection of random variables $\big\{\sqrt {\frac {n}{4\log n}}\Big(\frac 2 {J_A} - \EE\Big(\frac 2 {J_A}\Big)\Big)\big\}, 1\leq i \leq r, |A|=i,$ fulfils the assumptions of Proposition \ref{lemma_ind_thin}. First we make use of the main result of \cite{DK13} which says that the (centred and rescaled) internal lengths of different orders are asymptotically normal and independent. Together with \eqref{L_A_to_to_2_J_A_without_F} it follows that
\[
\Big(\sqrt {\frac {n}{4\log n}} \sum_{|A|=i}\Big(\frac 2 {J_A} - \EE\Big(\frac 2 {J_A}\Big)\Big)\Big)_{i=1, \dots, r} \stackrel d \longrightarrow N(0, I_r)
\]
where $I_r$ is the $r \times r$ identity matrix. The fact that the other two assumptions of Proposition \ref{lemma_ind_thin} hold has been proved in Lemma \ref{lemma_claims}. Therefore, we can apply Proposition \ref{lemma_ind_thin} for the collection of random variables  $\Big\{\sqrt {\frac {n}{4\log n}}\Big(\frac 2 {J_A} - \EE\Big(\frac 2 {J_A}\Big)\Big)\Big\}, 1\leq i \leq r, |A|=i,$ and the independent random coefficients $\{\alpha_1 V^A_{h_1}+ \alpha_2 V^A_{h_2}\}, 1\leq i \leq r, |A|=i$. 

Letting now $n \to \infty$ and using again Proposition \ref{proposition_bigger_r}  we obtain that the linear combination of the external lengths at the times $h_1$ and $h_2$ is asymptotically normal distributed. Note that we may interchange the limits in $r$ and $n$ since Proposition \ref{proposition_bigger_r} gives an uniform estimate in $r$. Therefore the stationary process $\ML^{n}$ converges in finite distributions to a stationary Gaussian process. It remains to compute its covariance function. 

Observe that for $h_1=0$ and $h_2=h$ 
\[\alpha_1V^A_{h_1}+\alpha_2V^A_{h_2}= \alpha_1\delta_{1,|A|}+ \alpha_2V^A_{h}\]
with second moment from \eqref{BD_prob_i} 
\begin{align*}
m_A =\begin{cases} 
      \alpha_1^2 +2\alpha_1\alpha_2 \Big(\frac 2 {h+2}\Big)^2 +\alpha^2_2\Big(\frac 2 {h+2}\Big)^2, &\text{ if } |A|=1 \\ 
			\alpha_2^2i\cdot\Big(\frac h {h+2}\Big)^{i-1}\Big(\frac 2 {h+2}\Big)^2 &\text{ if } |A|=i>1.
			\end{cases}
\end{align*}
Therefore from Proposition \ref{lemma_ind_thin} as $n \to \infty$
\begin{align*}
\sqrt{\frac {n}{4 \log n} }& \sum_{ |A| \leq r}\Big(\alpha_1 V^A_{0}+ \alpha_2 V^A_{ h}\Big) \cdot \Big(\frac 2 {J_A} - \EE\Big(\frac 2 {J_A}\Big)\Big)
\\
&\stackrel d \longrightarrow N\Big(0, \alpha_1^2 +2\alpha_1\alpha_2 \Big(\frac 2 {h+2}\Big)^2 + \alpha_2^2\Big(\frac 2 {h+2}\Big)^2\sum_{i=1}^r i\cdot\Big(\frac h {h+2}\Big)^{i-1}\Big).
\end{align*}
Taking the limit $r \to \infty$ we obtain that
\begin{align*}
\sqrt{\frac {n}{4 \log n} } &\Big(\alpha_1 (\ML_{0}^{n} - \EE(\ML_{0}^{n}))+ \alpha_2 (\ML_{h}^{n}- \EE(\ML_{h}^{n})\Big) \stackrel d \longrightarrow N\Big(0, \alpha_1^2 +2\alpha_1\alpha_2 \Big(\frac 2 {h+2}\Big)^2 + \alpha_2^2\Big) \qquad \text{as $n \to \infty$}.
\end{align*}
This gives the covariance function.

The almost sure continuity of the paths follows from Lemma 6.4.6 in \cite{MR06} using stationarity. To apply the lemma one has to check that there exists a $\delta>0$ such that 
\[\int_0^\delta \frac{\sigma^+(u)}{u \Big(\log \frac{1}{u}\Big)^{\nicefrac{1}{2}}}du < \infty,\]
where
\[\sigma^+(u):= \sup_{|t-s| \leq u \atop t,s \in [-\nicefrac{1}{2}, \nicefrac{1}{2}]} (\EE(\ML^1_t-\ML^1_s)^2)^{\nicefrac{1}{2}}.\]
This criterion is easily verified in our case, which finishes the proof.

\end{document}